\documentclass[12pt]{article}
\usepackage{amsfonts,amssymb,latexsym}
\newtheorem{theorem}{Theorem}[section]
\newtheorem{lemma}[theorem]{Lemma}
\newtheorem{corollary}[theorem]{Corollary}
\newtheorem{proposition}[theorem]{Proposition}

\newenvironment{proof}
{\par\addvspace{0.3cm}\noindent{\rm Proof. }}
{\nopagebreak\mbox{}\hfill $\Box$\par\addvspace{0.25cm}}

\newcommand{\be}{\begin{equation}}
\newcommand{\ee}{\end{equation}}

\newcommand{\bqn}{\begin{eqnarray}}
\newcommand{\eqn}{\end{eqnarray}}
\newcommand{\nn}{\nonumber}
\newcommand{\ba}{\begin{array}}
\newcommand{\ea}{\end{array}}
\newcommand{\tr}{\mathrm{trace}}
\renewcommand{\Re}{\mathrm{Re\,}}
\renewcommand{\Im}{\mathrm{Im\,}}
\newcommand{\wt}{\widetilde}

\newcommand{\iv}{^{-1}}
\newcommand{\iy}{\infty}
\newcommand{\ovl}{\overline}
\newcommand{\twomat}[1]{\left(\ba{cc} #1 \ea\right)}

\newcommand{\R}{{\mathbb R}}
\newcommand{\C}{{\mathbb C}}
\newcommand{\Z}{{\mathbb Z}}
\newcommand{\T}{{\mathbb T}}
\newcommand{\cM}{\mathcal{M}}
\newcommand{\cS}{\mathcal{S}}

\newcommand{\cC}{\mathcal{C}}
\newcommand{\W}{\mathcal{W}}

\newcommand{\cF}{\mathcal{F}}

\newcommand{\cU}{\mathcal{U}}
\newcommand{\cB}{\mathcal{B}}
\newcommand{\cL}{\mathcal{L}}

\newcommand{\ta}{\tilde{a}}
\newcommand{\tb}{\tilde{b}}
\newcommand{\tc}{\tilde{c}}
\newcommand{\Qh}{Q}
\newcommand{\uh}{\hat{u}}
\newcommand{\ah}{\hat{a}}

\newcommand{\Hh}{H_{\R}}

\begin{document}
\date{}
%\date{Jan 12, 2010}
\title{The asymptotics a Bessel-kernel determinant which arises in Random Matrix Theory
%\thanks{Keywords and phrases: random matrix, gap probability, determinant asymptotics, Wiener-Hopf operator, Hankel operator; MSC 2010: 47B35}
}
\author{Torsten Ehrhardt\thanks{Research supported in part by NSF Grant DMS-0901434.}\\  
Department of Mathematics \\ University of California \\ Santa Cruz, CA-95064, USA}
\maketitle

%%%

\begin{abstract}
In Random Matrix Theory the local correlations of the Laguerre and Jacobi Unitary Ensemble
in the hard edge scaling limit can be described in terms of the Bessel kernel
$$
B_\alpha(x,y) =\sqrt{x y}\, \frac{J_{\alpha}(x) y J'_\alpha(y)-
J_{\alpha}(y)x J'_\alpha(x)}{x^2-y^2},\qquad x,y>0,\quad \alpha>-1.
$$
In particular, the so-called hard edge gap probabilities $P^{(\alpha)}(R)$ can be expressed as the Fredholm determinants of the corresponding integral operator $B_\alpha$ restricted to the finite interval $[0,R]$. Using operator theoretic methods we are going to compute their asymptotics as $R\to\iy$, i.e., 
we show that 
$$
P^{(\alpha)}(R):=
\det(I-B_\alpha)|_{L^2[0,R]}
\sim
\exp\Big(-\frac{R^2}{4}+\alpha R-\frac{\alpha^2}{2}\log R\Big)
%e^{-\frac{R^2}{4}+\alpha R}R^{-\frac{\alpha^2}{2}}
\frac{G(1+\alpha)}{(2\pi)^{\alpha/2}},
$$
where $G$ stands for the Barnes $G$-function. In fact, this asymptotic formula will be proved for all complex parameters $\alpha$ satisfying $|\Re \alpha|<1$.
\end{abstract}

%%%%

\section{Introduction}
\label{s1}

For $1\le p\le \iy$ and a subinterval $I$ of the real line, let $L^p(I)$ stand for the
usual Lebesgue spaces. We denote by $B_\alpha$ the integral operator acting on $L^2(\R_+)$,
$\R_+=(0,\iy)$, with the kernel
\be\label{Bessel.kernel}
B_\alpha(x,y) =\sqrt{x y}\, \frac{J_{\alpha}(x) y J'_\alpha(y)-
J_{\alpha}(y)x J'_\alpha(x)}{x^2-y^2},
\ee
where $J_\alpha(x)$ are the Bessel functions with (complex) parameter $\alpha$. We will always assume that $\Re \alpha >-1$. 

In this paper we will be concerned with the Fredholm determinant of $B_\alpha$ on the interval $[0,R]$,
i.e., with the quantity
\bqn
P^{(\alpha)}(R) &:=& 1+\sum_{n=1}^\iy
\frac{(-1)^n}{n!}\int_0^R\dots\int_0^R\det\Big[B_\alpha(x_i,x_j)\Big]_{i,j=1}^n \, dx_1\dots dx_n.
\eqn
An equivalent definition can be given as an operator determinant \cite{GK,S},
\be\label{P.alpha}
P^{(\alpha)}(R) := \det(P_R-P_RB_{\alpha}P_R)|_{L^2[0,R]}.
\ee
Here $P_R$ stands for the orthogonal projection from $L^2(\R_+)$ onto the subspace $L^2[0,R]$,
\be\label{f.PR}
P_R:f(x)\mapsto g(x)=\left\{\ba{cl} f(x) & \mbox{ if } 0\le x\le R\\
0 & \mbox{ if } x>R.\ea\right.
\ee
We remark that the operators $P_R B_\alpha P_R|_{L^2[0,R]}$ are trace class operators on $L^2[0,R]$. 

The goal of this paper is to determine the asymptotics of
 $P^{(\alpha)}(R)$ as $R\to\iy$. In fact, we will prove the asymptotic formula 
\be\label{f.asym}
P^{(\alpha)}(R)
\sim
\exp\Big(-\frac{R^2}{4}+\alpha R-\frac{\alpha^2}{2}\log R\Big)
%e^{-\frac{R^2}{4}+\alpha R}R^{-\frac{\alpha^2}{2}}
\frac{G(1+\alpha)}{(2\pi)^{\alpha/2}},\qquad R\to\iy,
\ee
under the assumption $|\Re\alpha|<1$. Herein $G$ stands for the Barnes $G$-function \cite{Bar},
which is an entire function defined by
\be\label{Barnes}
G(1+z)=(2\pi)^{z/2}e^{-(z+1)z/2-\gamma_E z^2/2}
\prod_{k=1}^\iy \left((1+z/k)^{k} e^{-z+z^2/(2k)}\right)
\ee
with $\gamma_E$ being Euler's constant.

%%%%%%%%%

\vspace{4ex}

For real $\alpha>-1$, the quantity $P^{(\alpha)}(R)$ appears in Random Matrix Theory \cite{Me,TW,TW94b,TW-d}
as the gap probability for certain random matrix ensembles in the so-called hard edge scaling limit. Recall that the gap probability is the probability that no eigenvalue of the random matrix lies in an interval of some length.
The Bessel kernel $B_\alpha(x,y)$ arises, for instance, as the correlation function in the hard edge scaling limit of the Laguerre and Jacobi Unitary Ensemble (LUE/JUE) \cite{Fo93,NF95,NW93} as well as of generalized LUEs and JUEs \cite{KV02,V07}. All these ensembles consist of  complex self-adjoint matrices equipped with a certain probability measure which is invariant under unitary transform. 

To be more specific, recall that the LUE consists of positive self-adjoint complex $N\times N$ random matrices such that the joint probability density function of the eigenvalues is given by
$$
\rho_{\alpha,N}^{\mathrm{Lag}}(\lambda_1,\dots,\lambda_N)= c_{\alpha,N}\prod_{k=1}^N \lambda_k^{\alpha}e^{-\lambda_k} \prod_{1\le j<k\le N}|\lambda_j-\lambda_k|^2,\qquad \lambda_1,\ldots,\lambda_N>0.
$$
The JUE consists of all contractive self-adjoint complex $N\times N$ random matrices with joint probability density function of the eigenvalues given by
$$
\rho_{\alpha,\beta,N}^{\mathrm{Jac}}(\lambda_1,\dots,\lambda_N)= c_{\alpha,\beta,N}\prod_{k=1}^N (1-\lambda_k)^{\alpha}(1+\lambda_k)^{\beta} \prod_{1\le j<k\le N}|\lambda_j-\lambda_k|^2,\qquad -1<  \lambda_1,\ldots,\lambda_N<1.
$$
In both cases and for finite $N$, the probability that no eigenvalue lies in a subinterval $I$
of $\R_+$ or $[-1,1]$, respectively, can be written as a Fredholm determinant
$$
\det(I-K_N)|_{L^2(I)}= 1+\sum_{k=1}^\iy \frac{(-1)^k}{k!}
\int_{I}\dots\int_{I} \det \Big[K_N(x_r,x_s))\Big]_{1\le r,s\le k}dx_1\cdots dx_k.
$$
Here $K_N$ stands for the orthogonal projection onto the subspace of $L^2(\R_+)$ or $L^2[-1,1]$, resp., spanned by the first $N$ Laguerre or Jacobi functions, respectively, and $K_N(x,y)$ is the integral kernel of this operator.

In the random matrix context, $\alpha$  and $\beta$ are real parameters greater than $-1$. Let $P^{\mathrm{Lag},(\alpha)}_N(R)$ and $P^{\mathrm{Jac},(\alpha,\beta)}_N(R)$ denote the probabilities that  no eigenvalues lie in the interval $[0,R]\subset \R_+$ (Laguerre case)  or $[1-R,1]\subset[-1,1]$ (Jacobi case), respectively. Notice that these probabilities describe also the behavior of the eigenvalue closest to the hard edges of the ensembles.
With the appropriate scaling these probabilities converge (as $N\to\iy$) to the quantity $P^{(\alpha)}(R)$:
$$
P^{(\alpha)}(R)=\lim_{N\to\iy} P^{\mathrm{Lag},(\alpha)}_N\Big(\frac{R^2}{4N}\Big)=\lim_{N\to\iy} P^{\mathrm{Jac},(\alpha,\beta)}_N\Big(\frac{R^2}{2N^2}\Big)
$$
In fact, the Laguerre and Jacobi kernel converge, after the hard edge rescaling, to the Bessel kernel
$B_\alpha(x,y)$. This is also true for certain modified Laguerre and Jacobi random matrix ensembles \cite{KV02,V07}.

%%%%%%%%%%%%

\vspace{2ex}
Usually the kernel
$$
\widetilde{B}_\alpha(x,y)=\frac{J_{\alpha}(\sqrt{x}) \sqrt{y} J'_\alpha(\sqrt{y})-
J_{\alpha}(\sqrt{y})\sqrt{x} J'_\alpha(\sqrt{x})}{2(x-y)}
$$
is referred to as the Bessel kernel. It relates to our Bessel kernel by a simple change of variables.
In particular,
$$
P^{(\alpha)}(\sqrt{R})=\det(P_R- P_R \widetilde{B}_\alpha P_R)|_{L^2[0,R]}.
$$
Our form of Bessel kernel is more appropriate for the computations in this paper.

%%%%%%%
\pagebreak[1]
\vspace{4ex}

The asymptotics (\ref{f.asym}) was conjectured by Tracy and Widom \cite{TW94b}. We conjecture that it is valid for all complex $\alpha$ with $\Re\alpha>-1$, although our proof works only for $|\Re\alpha|<1$.
Following \cite{TW94b}, let us explain what has supported the Tracy-Widom conjecture. If we define
$$
\sigma(s)=-s\,\frac{d}{ds}\Big(\log P^{(\alpha)}(\sqrt{s})\Big),
$$
then it is known that $\sigma$ satisfies the differential equation
\be\label{f.PIII}
(s\sigma'')^2+\sigma'(\sigma- s\sigma')(4\sigma'-1)=\alpha^2(\sigma')^2,
\ee
which is reducible to a Painlev\'e III equation with specific parameters (see (3.13) in \cite{Jim}, and \cite{Oka}), and, in fact, also to a Painlev\'e V equation. (The equation (\ref{f.PIII}) is essentially the Jimbo-Miwa-Okamoto Painlev\'e III$'$ form.) Based on this differential equation, one can make a guess that the specific solution has
an asymptotic expansion as $s\to\iy$ in powers of $s^{1/2}$, i.e., 
$\sigma(s)=c_1s +c_2 s^{1/2}+c_3+c_4s^{-1/2}+\dots$. A straightforward computation yields recursively the coefficients (up to some ambiguity in a sign), and one obtains
\bqn\label{f.asexp}
\sigma(s)=\frac{s}{4}-\tau \frac{\alpha}{2} s^{1/2}+\frac{\alpha^2}{4}+\tau \frac{\alpha}{16} s^{-1/2}+\dots,\quad s\to\iy,
\eqn
with either $\tau=1$ or $\tau=-1$. The actual value of $\tau$ must be determined differently.
Finally, upon integration one obtains the asymptotics (\ref{f.asym}) for $P^{(\alpha)}(R)$, and even
higher order terms, {\em except for the constant term} $(2\pi)^{-\alpha/2}G(1+\alpha)$.

The conjecture \cite{TW94b} for the constant term (as well as for the correct value $\tau=1$ in (\ref{f.asexp})) relied on the special cases $\alpha=\pm1/2$. In these two cases, the Bessel operator equals a Wiener-Hopf-Hankel operator,
$$
B_{\pm 1/2}(x,y) = \frac{\sin(x-y)}{\pi(x-y)}\mp \frac{\sin(x+y)}{\pi(x+y)},
$$
and the asymptotics, including the constant, was conjecturally identified by Dyson \cite{Dy} based on rigorous work of Widom \cite{W1}. The rigorous identification of the constant (and in fact a proof of the asymptotics (\ref{f.asym})) in the cases $\alpha=\pm 1/2$ was given by the author \cite{E2} by employing
the same kind of methods that will be used in the present paper. Another proof for the special cases
$\alpha=\pm 1/2$ based on Painlev\'e transcendents and using the Riemann-Hilbert method \cite{DIZ}
was given recently by Baik, Buckingham, DiFranco, and Its \cite{BBdFI}. The quantities $P^{(\pm1/2)}(R)$ also occur in connection with the gap probability in the bulk scaling limit of the Gaussian Orthogonal and Symplectic Ensemble 
\cite{Me, TW-d}.
 
So far the derivation of (\ref{f.asexp}) is only heuristic, and the author is not aware of a rigorous proof. However, let us make the
following observation. Because it is comparatively easy to determine the asymptotics of $\sigma(s)$ as $s\to0$, the solution of the so-called connection problem for the Painlev\'e equation (\ref{f.PIII}) would provide such a proof. Unfortunately, the author did not find a rigorous solution of the connection problem for this concrete Painlev\'e equation. Thus this problem still seems to be open. However, considerable progress has been made in solving connection problems for other classes of Painlev\'e equations. The techniques employed there appear to be powerful enough to solve also the connection problem for (\ref{f.PIII}). We will not
discuss this topic any further, but refer the interested reader to the monograph of Fokas, Its, Kapaev, and Novokshenov \cite{FIKN}, as well as to the various articles of Kitaev, Andreev, Vartanyan, and Novokshenov referenced in this book.
 
There exist similar ``gap probability problems'', involving different kernels and arising from different scalings and/or random matrix ensembles. These problems were open for quite some time and were solved only recently. Besides the special case $\alpha=\pm1/2$ already mentioned above, the perhaps most important one was the problem of
asymptotics of the Fredholm determinant of the Wiener-Hopf sine kernel, $\sin(x-y)/(\pi(x-y))$. This quantity describes the gap probablity for the Gaussian Unitary ensemble in the bulk scaling limit.
The asymptotics was rigorously determined by Krasovsky \cite{Kr} using the Riemann Hilbert method and independently by the author in \cite{E1} using Wiener-Hopf and Toeplitz-Hankel operator techniques. Another proof was given by Deift, Its, Krasovsky, and Zhou \cite{DIKZ} also using the Riemann-Hilbert method. 

Another class of examples arises from the soft edge scaling of various random matrix ensembles.
Such scaling leads to the Airy kernel \cite{TW-a}. The corresponding asymptotics in the unitary case
was proved by Deift, Its, and Krasovsky \cite{DIK-a}, and another proof also including the orthogonal and
symplectic case was provided by Baik, Buckingham,  and DiFranco \cite{BBdF}. 

%%%%%%%%%%%%%%
%%%%%%%%%%%%%%

\vspace{4ex}

In order to give an outline of the paper and to make connections with other results, we need to
introduce some notation. The Fourier transform $\cF$ on $L^2(\R)$ and its inverse
will be written as 
\be\label{f.Fou}
(\cF f)(x) =\frac{1}{2\pi}\int_{-\iy}^\iy f(\xi) e^{-ix\xi}\, d\xi,
\qquad (\cF\iv g)(x)=\int_{-\iy}^\iy g(\xi)e^{ix \xi}\, d\xi.
\ee
The operators of primary interest to us are Wiener-Hopf and Hankel operators acting on
$L^2(\R_+)$. They are defined as
\bqn\label{f.WHH1}
W(a) &=& \Pi_+ \cF M_\R(a) \cF\iv  \Pi_+|_{L^2(\R_+)},
\\[1ex]
H_\R(a)  &=& \Pi_+\cF M_\R(a)\cF\iv  J_\R \Pi_+|_{L^2(\R_+)}.
\label{f.WHH2}
\eqn
Here $M_\R(a):f\mapsto a f$ is the multiplication operator on $L^2(\R)$ by a function $a\in L^\iy(\R)$,
$\Pi_+$ stands for the orthogonal projection of $L^2(\R)$ onto the subspace $L^2(\R_+)$,
and  $J_\R$ is the flip operator $(J_\R f)(x)=f(-x)$. Basic information about Wiener-Hopf and
Hankel integral operators can be found in \cite[Ch.~9]{BS}, but we mention that the notation is slightly
different there.

In addition to the finite section projection $P_R$ already defined in (\ref{f.PR}) we define
the complementary projection $Q_R=I-P_R$ acting on $L^2(\R_+)$.
Finally, for complex parameters $\beta$ we introduce the following two  functions defined on $\R$,
\be\label{f.uh}
\uh_{\beta,0}(x)=\left(\frac{ix+1}{ix-1}\right)^{\beta},\qquad
\uh_{\beta,\iy}(x)=\left(\frac{1+ix}{1-ix}\right)^{\beta},\qquad x\in\R.
\ee
Here the principal values of the power functions are considered. The function $\uh_{\beta,0}(x)$ is continuous on the one-point compactification  $\dot{\R}$ of $\R$ except at $x=0$, where it has a jump discontinuity. The function $\uh_{\beta,\iy}(x)$ is continuous on $\R$ and has (in general, different) limits
at $x=\pm \iy$.

%%%%%%%%
\vspace{4ex}

Let us now outline the proof of the asymptotic formula (\ref{f.asym}). It is split into several steps.
\vspace*{2ex}

{\em Step 1:} Here we will establish the following identity for each (fixed) $R>0$
 under the assumption $|\Re\alpha|<1$,
\bqn\label{ff.6}
P^{(\alpha)}(R)= \exp\Big(-\frac{R^2}{4}+\alpha R\Big)
\det\Big[\Big( P_R(I+H_\R(\psi))\iv P_R\Big)\iv\Big(
 P_R(I+H_\R(\hat{\psi}))\iv  P_R\Big)\Big]
\eqn
where
$$
\psi(x)=\uh_{-1/2-\alpha,0}(x)\uh_{1/2+\alpha,\iy}(x),\qquad
\hat{\psi}(x) =\uh_{-1/2,0}(x)\uh_{1/2+\alpha,\iy}(x)
$$
(see Cor.~\ref{c6.3}).
The expression under the determinant is of the form identity plus trace class operator on $L^2[0,R]$, and thus the determinant can be understood as an operator determinant \cite{GK, S}.
The several inverses appearing in this expression are those of bounded linear operators
acting on $L^2(\R_+)$ or $L^2[0,R]$, resp., whose existence will be shown.

The proof of identity (\ref{ff.6}) is carried out again in several steps (Sections \ref{s4}-\ref{s6}).
The main idea is to discretize the 
Bessel kernel by the Jacobi kernel (Sec.~\ref{s5}). In terms of random matrices, this corresponds precisely to taking the hard edge scaling limit of the JUE. In this way,
one obtains an interesting limit relation involving Hankel determinants,
$$
P^{(\alpha)}(R)=\exp\left(-\frac{R^2}{4}+\alpha R\right)\lim_{n\to\iy}
\frac{\det\left(\int_{-1}^1x^{j+k} \hat{w}_n(x)\, dx\right)_{j,k=0}^{n-1}}
{\det\left(\int_{-1}^1x^{j+k} \hat{w}(x)\, dx\right)_{j,k=0}^{n-1}},
$$
where $R>0$ is fixed,  and $\hat{w}_n$ and $\hat{w}$ are certain weight functions  on $[-1,1]$ depending also on $R$ and $\alpha$ (see Prop.~\ref{p5.2}). Remarkably, this identity 
exhibits already the leading order terms of the asymptotics of $P^{(\alpha)}(R)$ as $R\to\iy$.

Earlier, in Sec.~\ref{s4} we will establish a formula for Hankel determinants, which expresses them in terms of the discrete version of the types of operators appearing in  (\ref{ff.6}). Such formulas have been established before in special cases \cite{E1, E2}.

In Sec.~\ref{s6}, using these formulas, we then take the limit $n\to\iy$ while keeping $R>0$ fixed, in order to arrive at (\ref{ff.6}). We remark that the weight in the Hankel determinants has singularities at $-1$ and $1$, which is in a way the reason why we arrive at the rather complicated expression (\ref{ff.6}).

{\em Step 2:} In this step we now focus on the expression
\be\label{ff.6b}
\det\Big[\Big( P_R(I+H_\R(\psi))\iv P_R\Big)\iv\Big(
 P_R(I+H_\R(\hat{\psi}))\iv  P_R\Big)\Big].
\ee
One is tempted to write it as a quotient of two determinant. However, this is not possible because the underlying expressions are no longer of the form identity plus trace class. The reason is that both symbols $\psi$ and $\hat{\psi}$ are functions with a jump discontinuity at infinity. However, the ``size'' of the jumps is the same for $\psi$ and $\hat{\psi}$ so that in the expression (\ref{ff.6b}) some kind of cancellation occurs. 

What we will show is that (\ref{ff.6b})  equals asymptotically ($R\to\iy$)
$$
2^{-\alpha/2-\alpha^2}\,\frac{\det\Big( P_R(I+H_\R(\uh_{-1/2,0}))\iv P_R\Big)}{\det \Big(P_R(I+H_\R(\uh_{-1/2-\alpha,0}))\iv P_R\Big)}. 
$$
This result will be stated only the final section (Sec.~\ref{s8}), but the theorem leading to it is Thm.~\ref{t7.5}, which is proved in Sec.~\ref{s7}. The factor in front of the above quotient
can be interpreted as some correlation with the part $\uh_{1/2+\alpha,\iy}$ in the original symbols
$\psi$ and $\hat{\psi}$, which has been cancelled.

{\em Step 3:}
Here we apply a result, established by E.L.~Basor and the author \cite{BE3}, which computes the asymptotics of the determinants
$$
\det\Big( P_R(I+H_\R(\uh_{\gamma,0}))\iv P_R\Big)
$$
as $R\to\iy$. This quantity equals (up to a simple factor) a Wiener-Hopf-Hankel determinant with a
specific Fisher-Hartwig type symbol (see (\ref{ff.7}) below).
Notice that it is the asymptotics of this determinant which contributes the Barnes $G$-function  appearing in formula (\ref{f.asym}). With this last step the proof is complete  (Thm.~\ref{t8.2}).

In Step 1 and Step 2 we need two types of auxiliary results, which will be established in 
Sec.~2 and 3. One is the invertibility of certain operators, and the other one is the stability (i.e., asymptotic invertibility) of certain (generalized) sequences of operators.
These results will be applied at several different places, and they are certainly of
interest in their own rights. For this reason, we found it more suitable to establish these auxiliary results in two separate sections. Some notation and other basic auxiliary facts will also be stated in Sec.~2.

%%%%%%
\vspace{4ex}

Let us now make a connection with earlier results.
In \cite{E2} it has been shown that 
\bqn\label{ff.4}
P^{(-1/2)}(R)  &=& \exp\Big(-\frac{R^2}{4}-\frac{R}{2}\Big)\det\Big( P_R(I+H_\R(\hat{u}_{-1/2,0}))\iv P_R\Big),
\\[1ex] \label{ff.5}
P^{(1/2)}(R)  &=& \exp\Big(-\frac{R^2}{4}+\frac{R}{2}\Big)\det\Big( P_R(I-H_\R(\hat{u}_{1/2,0}))\iv P_R\Big).
\eqn
It is not too hard to see that both these identities are special cases of (\ref{ff.6}).

The types of determinants appearing in  (\ref{ff.4}) and (\ref{ff.5}), even for more general parameters,
can be identified up to a simple factor with Wiener-Hopf-Hankel determinants \cite{BE3}:
\be\label{ff.7}
\det   \Big( P_R(W(\hat{v}_{1/2+\alpha})+H_\R(\hat{v}_{1/2+\alpha})) P_R\Big)
= e^{-R(1/2+\alpha)}\det\Big(P_R(I+H_\R(\hat{u}_{-1/2-\alpha,0}))\iv P_R\Big)
\ee
if $|\Re\alpha|<-1$, and
\be\label{ff.8}
\det   \Big( P_R(W(\hat{v}_{-1/2+\alpha})-H_\R(\hat{v}_{-1/2+\alpha})) P_R\Big)
= e^{-R(-1/2+\alpha)}\det\Big(P_R(I-H_\R(\hat{u}_{1/2-\alpha,0}))\iv P_R\Big)
\ee
if $0<\Re\alpha<1$. Therein, $\hat{v}_\gamma(x)=(x^2/(1+x^2))^\gamma$. Notice that the right hand side of (\ref{ff.8}) makes sense for
$|\Re\alpha|<1$.

There is another way of looking at the determinant (\ref{P.alpha}) in the case of real $\alpha>-1$. In fact, the operator $I-B_\alpha$ can be considered as a Bessel convolution operator with highly
degenerate symbol. For $\alpha>-1$, the (unitary and self-adjoint) Hankel transform $\mathbf{H}_\alpha$ is defined by
$$
\mathbf{H}_\alpha:L^2(\R_+)\to L^2(\R_+), \quad
f(x)\mapsto g(x)=\int_0^\iy \sqrt{tx} J_\alpha(tx) f(t)\, dt.
$$
For $a\in L^\iy(\R_+)$ we define the Bessel convolution operator $B_\alpha(a)$ as
$$
B_\alpha(a)= \mathbf{H}_\alpha M_{\R_+}(a) \mathbf{H}_\alpha
$$
where $M_{\R_+}(a)$ stands for the multiplication operator on $L^2(\R_+)$ with symbol $a$.
If $a\in L^1(\R_+)\cap L^1(\R_+)$, then $B_\alpha(a)$ is an integral operator on $L^2(\R_+)$ with the
kernel
\be\label{Bessel.int}
B_\alpha(a)(x,y)=\int_0^\iy t\sqrt{xy} J_\alpha(tx)J_\alpha(ty)a(t)\, dt.
\ee

The Bessel operator $B_\alpha$ is a Bessel convolution operator where the symbol is
the characteristic function of the interval $[0,1]$, i.e., $B_\alpha=B_\alpha(\chi_{[0,1]})$. Indeed,
the kernel (\ref{Bessel.kernel}) can be rewritten as
\be\label{Bessel.int2}
B_\alpha(x,y)=\int_0^1  t\sqrt{xy} J_\alpha(tx)J_\alpha(ty)\, dt,
\ee
see, e.g., formula (25) in \cite{BE2b}. The quantities $P^{(\alpha)}(R)$ are thus the determinants of 
finite sections of the Bessel convolution operators with the symbol $1-\chi_{[0,1]}$, i.e., 
$$
P^{(\alpha)}(R)=\det  (P_{R} B_\alpha(1-\chi_{[0,1]}) P_R) |_{L^2[0,R]}.
$$
This symbol is highly singular in the sense that in vanishes on a whole interval. It is thus worse than Fisher-Hartwig type symbols and leads to a different asymptotics. This may be a reason why the asymptotics of the determinants is difficult to compute. The determinants of Bessel determinants with smooth and regular symbols have been computed in \cite{BE2b}, where the analogue of the Achiezer-Kac formula \cite{BS} was derived. Determinants of Bessel operators with Fisher-Hartwig type
symbols have not yet been investigated except for the cases $\alpha=\pm 1/2$ (see \cite{BE3,BEW}).

%%%%%%%%%%%%%%%%%%%%%%%%%%%%%%%%%%%%

\section{Notation and invertibility results}
\label{s2}
\subsection{Basic notation}
\label{s21}

For $1\le p\le \iy$ let $L^p(\T)$ stand for the Lebesgue spaces on the unit circle 
$\T=\{z\in\C\,:\, |z|=1\}$.  By $H^p(\T)$ and $\ovl{H^p(\T)}$  we denote the corresponding 
Hardy spaces, i.e.,
\bqn
H^p(\T) &=& \Big\{\; f\in L^p(\T)\;:\; f_n=0\mbox{ for all } n<0\;\Big\}, \\
\ovl{H^p(\T)} &=& \Big\{\; f\in L^p(\T)\;:\; f_n=0\mbox{ for all } n>0\;\Big\}.
\eqn
Here $f_n$ denotes the $n$-th Fourier coefficient of $f$, i.e, 
$$
f_n=\frac{1}{2\pi}\int_0^{2\pi} f(e^{ix}) e^{-in x}\, dx.
$$

Let us define the discrete analogues of the operators $W(a)$ and $H_\R(a)$.
For $a\in L^\iy(\T)$ the Toeplitz and Hankel operators acting on $H^2(\T)$ are defined by
\bqn\label{f.TH}
T(a)=P_\T M(a) P_\T|_{H^2(\T)},\qquad H(a)=P_\T M(a) J P_\T|_{H^2(\T)},
\eqn
where $P_\T$ is the orthogonal projection of $L^2(\T)$ onto the subspace $H^2(\T)$
(i.e., the Riesz projection), $J_\T$ is the flip operator $(J_\T f)(t)=t\iv f(t\iv)$, $t\in\T$, and
$M(a):f\mapsto af$ is the multiplication operator on $L^2(\T)$.
The following basic relations hold for $a,b\in L^\iy(\T)$,
\bqn
T(ab)= T(a) T(b) + H(a)H(\tilde{b}),\label{f.Tab}\\
H(ab)=T(a)H(b)+H(a) T(\tilde{b}).\label{f.Hab}
\eqn
Here and in what follows
\bqn
\tilde{b}(t):=b(t\iv),\quad t\in\T. 
\eqn
As special cases we obtain
\bqn\label{f.THabc}
T(abc)=T(a)T(b)T(c),\qquad H(ab\tilde{c})=T(a)H(b)T(c)
\eqn
for $a\in H^\iy(\T)$, $b\in L^\iy(\T)$, and $c\in \ovl{H^\iy(\T)}$.
For more information about Toeplitz and Hankel operators we refer to \cite[Ch.~2]{BS}.

We are going to consider Toeplitz and Hankel operators with particular symbols, which involve the  functions
\bqn\label{f.utau}
u_{\beta,\tau}(e^{i\theta} \tau) = e^{i\beta(\theta-\pi)},\qquad 0<\theta<2\pi, \quad \tau\in\T,\quad \beta\in\C,
\eqn
or, equivalently,
\bqn\label{f.utau2}
u_{\beta,\tau}(t)=(-t/\tau)^\beta,\qquad t\in\T
\eqn
with the principle value of the power function considered. These function are nonzero and continuous on $\T\setminus\{\tau\}$ and have a jump discontinuity at
$t=\tau$ with one-sided limits equal to $u_{\beta,\tau}(\tau\pm 0)=e^{\mp i \beta\pi}$. 

The finite section projections acting on $H^2(\T)$ are defined by
\bqn\label{f.Pn}
P_n:\sum_{k=0}^\iy {f_k}e^{ikx}\mapsto \sum_{k=0}^{n-1} f_k e^{ikx},
\eqn
and we put $Q_n=I-P_n$. There is a notational ambiguity with respect to the projections
$P_R$ and $Q_R$ acting on $L^2(\R_+)$, but the context will make clear what is meant.

The relationship between the discrete and the continuous case (except for the finite sections)
becomes clear when introducing the map $\cS: L^2(\T)\mapsto L^2(\R)$ defined as
the composition $\cS=\cF \circ \cU$, where $\cF$ is the Fourier transform and  $\cU$ is the unitary operator
$$
\cU: f\in L^2(\T)\mapsto g\in L^2(\R),\qquad g(x)=\frac{1}{\sqrt{\pi}(1-ix)}f\left(\frac{1+ix}{1-ix}\right).
$$
The restriction of $\cS$ onto $H^2(\T)$ is a mapping whose range can be identified with $L^2(\R_+)$. We will denote this restriction also by $\cS$. It is straightforward to show (see also \cite[Ch.~9]{BS})  that 
\be\label{HHR}
\cS T(a) \cS\iv = W(\ah),\quad
\cS H(a) \cS\iv = H_{\R}(\ah) \quad\mbox{with}\quad \ah(x)=a\left(\frac{1+ix}{1-ix}\right).
\ee
The transformation $a\mapsto \hat{a}$ also relates the specific symbols (\ref{f.uh})  and  (\ref{f.utau2}) 
to each other,
\be\label{uuh}
\uh_{\beta,0}(x)=u_{\beta,1}\left(\frac{1+ix}{1-ix}\right),\qquad
\uh_{\beta,\iy}(x)=u_{\beta,-1}\left(\frac{1+ix}{1-ix}\right).
\ee
We remark that up to a constant (due to our definition of the Fourier transform)
the operator $\cS$ is unitary.

%%%%%%%%%%%%%%%%%

\subsection{Invertibility of operators $I+H(\psi)$ and $I+H_\R(\hat{\psi})$}
\label{s22}

The goal of this section is to establish sufficient invertibility criteria for operators
$I+H(\psi)$ and $I+H_\R(\hat{\psi})$ acting on $H^2(\T)$ and $L^2(\R_+)$, respectively,  
for particular piecewise continious symbols.

Let $\W$ stand for the Wiener algebra on $\T$, i.e., for the set of all
$a\in L^\iy(\T)$ such that 
$$
\|a\|_{\W}:=\sum_{n=-\iy}^\iy |a_n| <\iy,
$$
where $a_n$ are the Fourier coefficients of $a$. Define
$$
\W_+=\W\cap H^\iy(\T), \qquad  \W_-=\W\cap \ovl{H^\iy(\T)}.
$$
The sets $\W$ and $\W_\pm$ are Banach algebras with unit element. For an arbitrary Banach algebra $B$ with unit element we denote by $GB$ the group of invertible elements
in $B$.
 
The following result about the invertibility of $I+H(\psi)$ generalizes previous results established in
Sec.~3.2 of \cite{BE3}, Sec.~4.1 of \cite{E1}, and Thm.~3.2 of \cite{E2}. The present result can be generalized further, but the proof is more complicated \cite{BEp}.

In the proof given below we will use the notions of  Fredholm operators and essential spectrum
\cite{GK}. Recall that a bounded linear operator $A$ acting on a space Hilbert space $H$ is called Fredholm if it has a closed range and if its kernel, $\ker A$, and its cokernel, $H/ \mathrm{im} A$, are both finite dimensional.
The Fredholm index is the difference between the dimensions of the kernel and cokernel.
The essential spectrum of $A$ is the set of all complex numbers $\lambda$ such that 
$A-\lambda I$ is not a Fredholm operator.

\begin{theorem}\label{t21}
Let $\alpha,\beta\in \C$ be such that $|\Re\alpha|<1$ and $|\Re\beta|<1$, and let 
$c_+\in G\W_+$. Then the operator
$$
I+H(\psi)\quad \mbox{with}\quad \psi=\tc_+ c_+\iv  u_{-1/2-\alpha,1}u_{1/2+\beta,-1},
$$
is invertible on $H^2(\T)$.
\end{theorem}
\begin{proof}
We first use a result of Power \cite{Po1, Po2} (see also \cite[Sec.~4.7]{BS}) in order to show that the operator $I+H(\psi)$ is a Fredholm operator. The result of Power states that the essential spectrum of a Hankel operator $H(\psi)$ with piecewise continuous
symbol $\psi$ equals the following union of closed intervals in the complex plane,
$$
\mathrm{sp}_{\mathrm{ess}} H(\psi)=[0,i\psi_{-1}]\cup [0,-i \psi_1]\cup\bigcup\limits_{\tau\in \T \atop \Im(\tau)>0}\left[-i\sqrt{\psi_\tau\psi_{\bar{\tau}}}, i\sqrt{\psi_\tau\psi_{\bar{\tau}}}\,\right],
$$
where
$$
\psi_\tau=\frac{1}{2}(\psi(\tau+0)-\psi(\tau-0))\quad \mbox{ with }\quad \psi(\tau\pm0)=\lim_{\varepsilon\to+0} \psi(\tau e^{i\varepsilon}).
$$
In our case we have $\psi_\tau=0$ for $\tau\in \T$, $\Im(\tau)>0$, and
$$
\psi_1= i \cos(\alpha\pi),\qquad \psi_{-1}=-i \cos(\beta\pi).
$$
Hence the essential spectrum of $H(\psi)$ is $[0,\cos(\alpha\pi)]\cup[0,\cos(\beta\pi)]$.
The point $\lambda=-1$ does not  belong to the essential spectrum if and only if $\Re\alpha$ and $\Re\beta$ do not belong to the set $1+2\Z$. This is fulfilled in our case. Hence we conclude that $I+H(\psi)$ is Fredholm. Moreover, making use of the fact that the complement of the essential spectrum is connected and the fact that the Fredholm index is invariant under small perturbations it follows that the Fredholm index of $I+H(\psi)$ is zero.

It remains to show that the kernel of $I+H(\psi)$ is trivial. For $\tau\in\T$, we introduce the functions
$$
\eta_{\gamma,\tau}(t)=(1-t/\tau)^\gamma,\qquad \xi_{\gamma,\tau}(t)=(1-\tau/t)^{\gamma},\qquad t\in\T\setminus\{\tau\},
$$
where the principal values of the power functions are considered. We notice that 
$$
u_{\gamma,\tau}(t)=\xi_{-\gamma,\tau}(t)\cdot \eta_{\gamma,\tau}(t) \quad \mbox{and}\quad
{\eta}_{\gamma,\tau}(t\iv)=\xi_{\gamma,\tau\iv}(t).
$$
Now assume that $f_+\in H^2(\T)$ belongs to the kernel of $I+H(\psi)$. Then there exists $f_-\in \overline{H^2(\T)}$
such that 
$$
f_+(t)+\psi(t)t\iv \tilde{f}_+(t)=t\iv f_-(t).
$$
We rewrite this as
$$
t f_+(t)+\psi(t) \tilde{f}_+(t)= f_-(t),
$$
and decompose (using $u_{1,1}(t)=-t$)
$$
-t\psi=  \tc_+ c_+\iv u_{1/2-\alpha,1} u_{1/2+\beta,-1}
=
(c_+\iv \eta_{1/2-\alpha,1}\eta_{1/2+\beta,-1})\cdot
(\tc_+\iv  \xi_{1/2-\alpha,1}\xi_{1/2+\beta,-1}) \iv
$$ 
in order to obtain
$$
f_0:=\tc_+\iv \xi_{1/2-\alpha,1}\xi_{1/2+\beta,-1} t f_+-
c_+\iv \eta_{1/2-\alpha,1}\eta_{1/2+\beta,-1}t\iv \tilde{f}_+
=\tc_+\iv \xi_{1/2-\alpha,1}\xi_{1/2+\beta,-1} f_-.
$$
Because of the assumptions on $\alpha$ and $\beta$, the function $\tilde{c}_+\iv \xi_{1/2-\alpha,1}\xi_{1/2+\beta,-1}$ belongs to $\overline{H^2(\T)}$. Hence the right hand side of the above equation belongs to $\overline{H^1(\T)}$. Each of the terms on the left hand side belongs to $L^1(\T)$ and 
$\tilde{f}_0=-f_0$. Comparing the Fourier coefficients it follows that $f_0=0$, whence $f_-=0$ and
$$
f_+(t)=-t\iv\psi(t)\tilde{f}_+(t).
$$
Now we decompose (using $u_{-1,-1}(t)=t\iv$)
$$
t\iv \psi=
 \tc_+ c_+\iv u_{-1/2-\alpha,1} u_{-1/2+\beta,-1}
=
(c_+ \eta_{1/2+\alpha,1}\eta_{1/2-\beta,-1})\iv \cdot
(\tc_+  \xi_{1/2+\alpha,1}\xi_{1/2-\beta,-1})$$
in order to obtain
$$
g_+:= c_+ \eta_{1/2+\alpha,1}\eta_{1/2-\beta,-1}f_+ 
=-\tc_+  \xi_{1/2+\alpha,1}\xi_{1/2-\beta,-1} \tilde{f}_+.
$$
The function $g_+$ belongs to $H^1(\T)$, and the last equation states that $g_+=-\tilde{g}_+$.
It follows that $g_+=0$, and from this that $f_+=0$. Hence the kernel of $I+H(\psi)$ is trivial.
Since we have already shown that $I+H(\psi)$ is  a Fredholm operator with index zero, this
implies the invertibility.
\end{proof}

The analogue of the previous result in the continuous case is stated next. For simplicity, we restrict ourselves to symbols in which $c_+\equiv 1$,
which is enough for our purposes.

\begin{corollary}\label{c22}
Let $\alpha,\beta\in \C$ be such that $|\Re\alpha|<1$ and $|\Re\beta|<1$. 
Then the operator
$$
I+\Hh(\uh_{-1/2-\alpha,0} \uh_{1/2+\beta,\infty})
$$
is invertible in $L^2(\R_+)$.
\end{corollary}
\begin{proof}
We apply the transform $\cS:H^2(\T)\to L^2(\R_+)$ to the Hankel operator $H(\psi)$ and obtain 
$$
\cS H(u_{-1/2-\alpha,1}u_{1/2+\beta,-1})\cS\iv= H_{\R}(\uh_{-1/2-\alpha,0} \uh_{1/2+\beta,\iy})
$$
by using the formulas (\ref{HHR}) and (\ref{uuh}). Now the result follows from the previous theorem.
\end{proof}

%%%%%%%%%%%%%%%%%%%%%%%%%%%%%%%%%%%%

%%%%%%%%%%%%%%%%%%%%%%%%%%%%%%%%%%%%
%%%%%%%%%%%%%%%%%%%%%%%%%%%%%%%%%%%%
%%%%%%%%%%%%%%%%%%%%%%%%%%%%%%%%%%%%
 
\section{Several results about stability}
\label{s3}

\subsection{Definitions and basic results}
\label{s31}

We  need the notions of stability and strong convergence. These notions involve generalized sequences of operators. Let $\Lambda\subset\R$ be an index set and
assume that the supremum $\lambda_\iy=\sup\Lambda\in \R\cup\{+\iy\}$ does not belong to $\Lambda$.
For us, only the following three settings are of interest,
\begin{enumerate}
\item[(i)] 
$\Lambda=\Z_+$ and $\lambda_\iy=+\iy$,
\item[(ii)] 
$\Lambda=(0,\iy)$ and $\lambda_\iy=+\iy$,
\item[(iii)]
$\Lambda=[0,1)$ and $\lambda_\iy=1$.
\end{enumerate}

We consider generalized sequences of bounded linear operators 
$\{A_\lambda\}_{\lambda\in\Lambda}$, where $A_\lambda\in \cL(H_\lambda)$ and $H_\lambda$ are Hilbert spaces. Such a sequence is called {\em stable} if there exists a $\lambda_0\in\Lambda$ such that the operators $A_\lambda$ are invertible for all $\lambda\ge \lambda_0$, $\lambda\in \Lambda$,
and if the inverses are (uniformly) bounded, i.e., 
$$
\sup_{\lambda\ge \lambda_0}\|A_\lambda\iv\|_{\cL(H_\lambda)}<\iy.
$$
Stability is also referred to as asymptotic invertibility. 

Given a sequence $\{A_\lambda\}_{\lambda\in\Lambda}$, let us assume that the 
spaces $H_\lambda$ are subspaces of a possibly larger Hilbert space $H$, and let 
$A\in\cL(H)$. We say that $A_\lambda\to A$ converges {\em strongly} on $H$ as $\lambda\to\lambda_\iy$ if for each $x\in H$ we have
$$
\lim_{\lambda\to\lambda_\iy} \|A_\lambda P_\lambda x-Ax\|_{H}=0.
$$
Therein, $P_\lambda$ stands for the orthogonal projection from $H$ onto $H_\lambda$.

We will always be concerned with bounded sequences, i.e., sequences $\{A_\lambda\}_{\lambda\in\Lambda}$ for which
$$
\sup_{\lambda\in\Lambda}\|A_\lambda\|_{\cL(H_\lambda)}<\iy.
$$
Such a sequence is called  a zero sequence if
$$
\lim_{\lambda\to\lambda_\iy}\|A_\lambda\|_{\cL(H_\lambda)}=0.
$$
The above definitions and the following basic results can be found, e.g.,  in \cite[Sec.~7.1]{BS}.
The adjoint of a linear bounded operator $A$ acting on a Hilbert space will be denoted by
$A^*$.

\begin{lemma}\label{l.basic}
Let  $\{A_\lambda\}$ be a bounded sequence of operators $A_\lambda\in \cL(H_\lambda)$. Then
the following holds:
\begin{enumerate}
\item[(i)]  $\{A_\lambda\}$
is stable if and only if there exists a bounded sequence $\{B_\lambda\}$
such that $\{A_\lambda B_\lambda-I\}$ and $\{B_\lambda A_\lambda-I\}$ are zero sequences.
\item[(ii)]
If $\{A_\lambda\}$ is stable, and $A_\lambda\to A$ and $A_\lambda^*\to A^*$ strongly, then
$A$ is invertible.
\item[(iii)]
If  $\{A_\lambda\}$ is stable, $A$ is invertible,  and $A_\lambda\to A$ strongly,  then  
$A_\lambda\iv\to A\iv$ strongly.
\end{enumerate}
\end{lemma}

The significance of the previous lemma is the following.
Statement (i) implies that stability is invariant under perturbation by zero sequences. It also justifies to refer to stability as asymptotic invertibility. Statement (iii) indicates the reason why we are interested in stability. In fact, we frequently need to show the strong convergence of a sequence of the inverses of operators.
Finally, statement (ii) shows that in the cases of interest to us (in which we will always have strong convergence $A_\lambda\to A$ and $A_\lambda^*\to A^*$) the assumption in (iii) that $A$ be invertible is also necessary.

Since we will be dealing with operator determinants we need the notion of trace class operators \cite{GK,S}. Recall that a compact operator $A$ acting on a Hilbert space is called trace class
if the sequence of its singular values (i.e., the eigenvalues of $(A^*A)^{1/2}$) is absolutely summable.
The sum of the singular values is by definition the trace norm of $A$.

We often want to show convergence in the trace norm. Here the following
basic statement is the key (see also \cite[Sec.~1.3]{BS}).

\begin{lemma}\label{l.strongconv}
Assume that $K$ is a trace class operator and that $A_\lambda\to A$ and $B_\lambda^*\to B^*$ strongly. Then $A_\lambda K B_\lambda \to AKB$ in the trace norm.
\end{lemma}

In order to show the strong convergence of certain operator sequences we will frequently employ the following result (see, e.g., \cite[Lemma 5.4]{ES1}).
Therein, the measure is meant to be the Lebesgue measure.

\begin{lemma}\label{l3.3}
Let $K=\R$ or $K=\T$.
Let $\{b_\lambda\}_{\lambda\in \Lambda}$ be a bounded sequence of functions in $L^\iy(K)$, and assume that $b_\lambda$ converges in measure to $b\in L^\iy(K)$. Then
$$
W(b_\lambda)\to W(b),\quad H_\R(b_\lambda)\to H_\R(b),\quad
W(b_\lambda)^*\to W(b)^*,\quad H_\R(b_\lambda)^*\to H_\R(b)^*,
$$
strongly on $L^2(\R_+)$ in the case $K=\R$, and,
$$
T(b_\lambda)\to T(b),\quad H(b_\lambda)\to H(b),\quad
T(b_\lambda)^*\to T(b)^*,\quad H(b_\lambda)^*\to H(b)^*,
$$
strongly on $H^2(\T)$ in the case $K=\T$, respectively.
\end{lemma}

In the following two subsections we are going to use the following auxiliary result as well (see \cite[Prop.~7.15]{BS}).

\begin{lemma}\label{l21.1}
Let $A$ be an invertible bounded linear operator on a Hilbert space $H$. Let $P$ be a bounded
projection operator on $H$, and let $Q=I-P$. Then $PAP$ is invertible on $\mathrm{im}\, P$ if and only if
$QA\iv Q$ is invertible on $\mathrm{im}\, Q$. In particular,
\bqn
(PAP)\iv  &=& PA\iv P- PA\iv Q (QA\iv Q)\iv QA\iv P,\nn\\
(QA\iv Q)\iv & =& Q A Q - Q A P (PAP)\iv P A Q\nn.
\eqn
\end{lemma}
To be more specific, in the previous lemmas the operators $PAP$ and $QA\iv Q$ are considered as
the restrictions onto the spaces $\mathrm{im}\, P$ and  $\mathrm{im}\, Q$, respectively.
 
 %%%%%%%%%%%%%%%%

\subsection{Stability of certain finite sections (continuous case)}
\label{s32a}

The goal of this subsection is to prove the stability of the following sequence of finite sections,
\bqn\label{seq1}
\{P_R (I+H_\R(\hat{\psi}))\iv P_R\}_{R>0},
\eqn
for certain symbols $\hat{\psi}\in PC(\R)$. Here $PC(\R)$ stands for the set of all piecewise continuous functions on the real line (with limits as $x\to\pm\iy$).

We remark that the main result of this subsection (Corollary \ref{c3.7}) could also be derived from
results of Roch, Santos, and Silbermann \cite{RSS}. However, some additional effort would still be necessary in order to identify certain operators which are defined by homomorphisms and to show the invertibility of these operators. In fact, these invertibility results would correspond to Corollary \ref{c22} and Proposition \ref{p22.1} below. We found it easier to give a direct proof in the special case that is of interest to us. 
Nonetheless, Corollary \ref{c3.7} can probably be generalized to a larger class of symbols, in which case a direct application of the results  \cite{RSS} is perhaps more suitable.

Introduce the function
\be\label{f.chi}
\chi(x)=i \cdot \mathrm{sign}(x)=\left\{\ba{rl} i &\mbox{ if } x>0\\ -i &\mbox{ if } x<0.\ea\right.
\ee
In a first step we establish an invertibility or, more precisely, spectral result.

\begin{proposition}\label{p22.1}
For  $R>0$, the operators $\Qh_R\Hh(\chi)\Qh_R|_{L^2{(R,\iy)}}$ are self-adjoint, unitarily equivalent to each other,  and have spectrum equal to the interval $[0,1]$.
\end{proposition}
\begin{proof}
Let us consider the Fourier convolution operator
$$
W_0(a)=\cF M_\R(a) \cF\iv,
$$
where we use the same notation as in (\ref{f.WHH2}). Then the operators $Q_R\Hh(\chi)Q_R$ are equal to the compressions of the operator
$W_0(\chi)  J_\R$ onto $L^2(R,\iy)$. It is well known (see, e.g., \cite{HRS})
that the Fourier convolution operator $W_0(\chi)$ is equal to $iS$, where $S$ is the singular integral operator on $\R$, 
$$
(S f)(x) = \frac{1}{\pi i } \int_{-\iy}^\iy \frac{f(y)}{x-y}\, dy, \qquad f\in L^2(\R).
$$
Therein the integral exists a.e.\ as the Cauchy principle value. It follows that 
the operator $\Hh(\chi)$ is an  integral operator on $L^2(\R_+)$ with the integral kernel
$\pi \iv (x+y)\iv$. 

On the other hand, $H_\R(\chi)$ is a Mellin convolution operator $M_0(n)$ with the generating function
$$
n(z)=
(\cosh(\pi z))\iv,\qquad z\in\R.
$$ 
To see this, recall that a Mellin convolution operator is defined as $M_0(a)=\cM\iv M_\R(a) \cM$, where
$$
\cM:L^2(\R_+)\to L^2(\R),\qquad (\cM f)(\xi)= \int_{0}^\iy f(x) x^{i\xi-1/2}\, dx
$$
is the Mellin transform. Obviously, $\mathcal{M}=\mathcal{F}\iv \mathcal{T}$, where 
$\mathcal{T}$ is the isometry defined by
$$
\mathcal{T}:L^2(\R_+)\to L^2(\R),\qquad
(\mathcal{T} f)(x)=f(e^x)e^{x/2}.
$$
Therefore the identity $H_\R(\chi)=M_0(n)$ reduces to
$$
\mathcal{T} H_\R(\chi)\mathcal{T}\iv = \cF M_\R(a) \cF\iv = W_0(n).
$$
This identity can be verified directly by showing that the Fourier transform of $n(z)$ is 
$(2\pi \cosh(x/2))\iv$ and by observing that $\mathcal{T} H_\R(\chi)\mathcal{T}\iv$
is an integral operator on $L^2(\R)$ with the integral kernel equal to
$(2\pi  \cosh((x-y)/2))\iv$.

Finally, observe that the transform $A\mapsto \mathcal{T}A\mathcal{T}\iv$ maps
the operator $Q_R$ into $M_\R(\chi_{(\ln R,\iy)})\in \cL(L^2(\R))$. Hence
the operators $Q_R H_\R(\chi) Q_R$ are unitarily equivalent to 
$$
M_\R(\chi_{(\ln R,\iy)})W_0(n)M_\R(\chi_{(\ln R,\iy)})|_{L^2(\ln R,\iy)},
$$
which (by means of a translation) are unitarily equivalent to the Wiener-Hopf
operator $W(n)$. In particular, they are unitarily equivalent to each other.
Because the symbol $n(z)$ is real-valued, continuous on the one-point compactification $\dot{\R}$, and has range $[0,1]$, it follows from basic Wiener-Hopf theory that the operator $W(n)$ is self-adjoint and has spectrum equal to $[0,1]$.
\end{proof}

The fact that  the operators $Q_RH_\R(\chi)Q_R|_{L^2(R,\iy)}$ are unitarily equivalent to each other
can also be seen in a more direct way. Introducing the unitary operators
$\hat{Y}_R:f(x)\mapsto R^{-1/2}f(x/R)$, acting on $L^2(\R_+)$, one can show easily that
\bqn\label{f.YRa}
Q_1 H_\R(\chi) Q_1= \hat{Y}_R\iv Q_R H_\R(\chi) Q_R \hat{Y}_R.
\eqn
Moreover, we can immediately conclude the following result.

\begin{corollary}\label{c3.4}
Let $\alpha\in\C$ be such that $|\Re\alpha|<1$. Then the operators
\bqn\label{f.31x}
P_R(I+H_\R(\hat{u}_{-1/2-\alpha,0}\hat{u}_{1/2+\alpha,\iy}))\iv P_R|_{L^2[0,R]}
\eqn
are invertible for each $R>0$ and the inverses are uniformly bounded.
\end{corollary}
\begin{proof}
The invertibility of $I+H_\R(\hat{u}_{-1/2-\alpha,0}\hat{u}_{1/2+\alpha,\iy})$ follows from 
Corollary \ref{c22}. Hence, using Lemma \ref{l21.1}, the invertibility of the above operators
is equivalent to the invertibility of
$$
Q_R(I+H_\R(\hat{u}_{-1/2-\alpha,0}\hat{u}_{1/2+\alpha,\iy}))Q_R.
$$
Now observe that 
$$
\hat{u}_{-1/2-\alpha,0}(x)\hat{u}_{1/2+\alpha,\iy}(x)=
\left\{\ba{cl} e^{i\pi(1/2+\alpha)} & \mbox{ if } x>0\\
e^{-i\pi(1/2+\alpha)} & \mbox{ if } x<0,
\ea\right.
$$
which is a constant function plus the function $\cos(\pi \alpha)\chi(x)$. The last operator thus equals
$$
Q_R+\cos(\pi\alpha) Q_R H_\R(\chi)Q_R,
$$
which is invertible by the previous proposition if and only if $\cos(\pi \alpha)\notin(-\iy,-1]$, i.e.,
$\Re\alpha\notin1+2\Z$. Since for different $R$ these operators are unitarily equivalent to each other,
the uniform boundedness of the inverses follows immediately.
Applying the formulas for the inverses stated in Lemma \ref{l21.1} 
it follows that the inverses of (\ref{f.31x}) are uniformly bounded, too.
\end{proof}

Let $\overline{\R}$ stand for the two-point compactification of $\R$, and let $C(\ovl{\R})$ stand for the set of all continuous functions on $\overline{\R}$.

\begin{lemma}\label{l22.2}
Let $a\in C(\overline{\R})$. Then $\|\Qh_R \Hh(a)\Qh_R\|_{\cL(L^2(R,\iy))}
\to 0$ as $R\to\iy$.
\end{lemma}
\begin{proof}
A Hankel operator $\Hh(a)$ with $a\in C(\dot{\R})$ is compact. Because $Q_R\to 0$ strongly
as $R\to\iy$, the statement is proved for such symbols. To extend the result to 
arbitrary functions in $C(\overline{\R})$ it suffices to consider the symbol 
$a(x)=\chi(x)-\chi(x)e^{-|x|}$, which is continuous on $\ovl{\R}$, but has different limits as $x\to\pm\iy$.
As shown in the proof of the Proposition \ref{p22.1}, $H_\R(\chi)$ is an integral operator on $L^2(\R_+)$ with integral kernel
$\pi\iv(x+y)\iv$. By computing the Fourier transform of $\chi e^{-|x|}$ we see that  $H_\R(\chi e^{-|x|})$
has the integral kernel $\pi\iv(1+x+y)\iv$. Hence $H_\R(a)$ has the kernel $(\pi (x+y)(1+x+y))\iv$. It is now easily seen that 
$\Qh_R\Hh(a)\Qh_R|_{L^2(R,\iy)}$ are Hilbert-Schmidt operators with the Hilbert-Schmidt norm converging to
zero as $R\to\iy$. This implies the desired assertion.
\end{proof}
 
The key point to the treatment of the stability of (\ref{seq1}) is the 
following result. Therein the generating functions are allowed to have jump discontinuities at
$x=0$ and $x=\iy$. Notice that the size of the jump at $x=\iy$ does  not play a role.

\begin{theorem}\label{t22.3}
Let $a\in PC(\R)$ be continuous on $\R\setminus\{0\}$. Then the generalized sequence 
$\left\{\Qh_R+\Qh_R H_\R(a) \Qh_R\right\}_{R>0}$ is stable if and only if
$$\beta:=\frac{a(+0)-a(-0)}{2i} \notin (-\iy,-1].$$ 
\end{theorem}
\begin{proof}
Given $a$ we write
$$
a(x)= \beta\,  \chi(x) + b(x)
$$
with $b\in C(\overline{\R})$. Because of the previous lemma, the sequence $\Qh_R H_\R(b)\Qh_R$ converges to zero. Hence the sequence under consideration is stable if and only if so is
$Q_R+\beta Q_R H_\R(\chi) Q_R$ (see Lemma \ref{l.basic}(i)). However, the necessary and sufficient condition for the stability of this sequence was already identified in the proof of Corollary \ref{c3.4} by help of Proposition 
\ref{p22.1}.
\end{proof}

The desired stability result is the following. We restrict ourselves to the cases which are of interest to us.

\begin{corollary}\label{c3.7}
Let $\alpha\in \C$ be such  that $|\Re \alpha|<1$. Then the generalized sequences
$$
\left\{ P_R(I+H_\R(\uh_{-1/2-\alpha,0}))\iv  P_R\right\}_{R>0},\qquad
\left\{ P_R(I+H_\R(\uh_{1/2+\alpha,\iy}))\iv  P_R \right\}_{R>0}
$$
are stable.
\end{corollary}
\begin{proof}
By Corollary \ref{c22}, the assumption on $\alpha$ implies that the operators $I+H_\R(\uh_{-1/2-\alpha,0})$
and $I+H_\R(\uh_{1/2+\alpha,\iy})$ are invertible on $L^2(\R_+)$. Hence we can apply Lemma \ref{l21.1} (including the formulas for the inversess) and reduce the stability statements to the stability of 
$$
\left\{\Qh_R(I+H_\R(\uh_{-1/2-\alpha,0}))\Qh_R\right\}_{R>0}\quad\mbox{and}\quad
\left\{\Qh_R(I+H_\R(\uh_{1/2+\alpha,\iy}))\Qh_R\right\}_{R>0},
$$
respectively. By Lemma \ref{l22.2} the second sequence is always stable, regardless of the parameter $\alpha$. For the first sequence we can apply the previous theorem, and therefore we compute the parameter
$$
\beta=
\frac{\uh_{-1/2-\alpha,0}(+0)-\uh_{-1/2-\alpha,0}(-0)}{2i}=\frac{e^{\pi i(1/2+\alpha)}-e^{-\pi i (1/2+\alpha)}}{2 i}=\cos(\pi\alpha).
$$
{}The condition $\cos(\alpha\pi)\notin(-\iy,-1]$, i.e., $\Re\alpha\notin 1+ 2\Z$, is satisfied by assumption.
\end{proof}

%%%%%%%%%%

\subsection{Stability of certain finite sections (discrete case)}
\label{s32b}

Now we are going to present the discrete versions of the previous results and prove the stability
of sequences
\bqn\label{seq2}
\{P_n(I+H(\psi)) \iv P_n\}_{n\in\Z_+}
\eqn
for certain symbols $\psi\in PC(\T)$. Here $PC(\T)$ stands for the
set of all piecewise continuous functions on $\T$.

Also in this case there exists a general stability criterion, which was established by Roch \cite{R87} and
from which the main result of this subsection (Cor.~\ref{c3.8}) could be derived. As before
additional effort would still be necessary. Therefore we found it easier to give a direct proof
for the particular cases we are interested in. A generalization of Corollary \ref{c3.8} to more general symbols is certainly possible.

\begin{proposition}\label{p3.8}
Let $a\in PC(\T)$ be continuous on $\T\setminus\{1\}$.
Then the sequence $\{Q_n+Q_nH(a)Q_n\}_{n\in\Z_+}$ is stable if and only if
$$\beta:=\frac{a(1+0)-a(1-0)}{2i} \notin (-\iy,-1].$$ 
\end{proposition}
\begin{proof}
Let $\sigma(e^{ix})=i(\pi-x)/\pi$, $0<x<2\pi$. Then we can write
$$
a(e^{ix}) =\beta   \sigma(e^{ix})+b(e^{ix})
$$
with $b\in C(\T)$. Because $H(b)$ is compact and $Q_n\to 0$ strongly,
the stability of $Q_n+Q_nH(a) Q_n$ is equivalent to the stability of
$Q_n+\beta Q_n H(\sigma) Q_n$ (see Lemma \ref{l.basic}(i)).

Let us introduce the operators
$$
E_n: H^2(\T)\to L^2(\R_+), \quad \sum_{k=0}^\iy f_k t^n \mapsto \sqrt{n}\sum_{k=0}^\iy f_k \chi_{[\frac{k}{n},\frac{k+1}{n}]}(x),
$$
which are isometries. Their adjoints are given by
$$
E_n^*: L^2(\R_+)\to H^2(\T), \quad f(x)\mapsto \sum_{k=0}^\iy f_k t^n\quad\mbox{with}\quad
f_k=\sqrt{n} \int_{k/n}^{(k+1)/n} f(x)\, dx.
$$

Now consider the operator $E_{1}^* H_\R(\chi) E_1$.  In the matrix representation of this operator with respect to the standard basis $\{t^n\}_{n=0}^\iy$ of $H^2(\T)$, the $(j,k)$-entry equals
$$
\frac{1}{\pi}\int_{j}^{j+1}\int_{k}^{k+1}\frac{1}{x+y}\, dx\, dy=\frac{1}{\pi (1+j+k)}+O((j+k+1)^{-2}), \qquad j,k\ge0.
$$
The error term corresponds to a Hilbert-Schmidt operator $K$. Because the $n$-th Fourier coeffcient
of $\sigma$ evaluates to $(\pi n)\iv$, $n\neq0$, we obtain from the definition of the Hankel operators
that  $E_1^* H_\R(\chi)E_1= H(\sigma)+K$.
From this we now conclude that  the stability of $Q_n + \beta Q_n H(\sigma) Q_n$ is equivalent to the stability of $Q_n + \beta Q_n E_1^* H_\R(\chi) E_1 Q_n$ (see Lemma \ref{l.basic}(i) and use the fact that $Q_n K\to 0$ as $n\to\iy$).

Now we observe that 
$$
Q_n E_1^* \Hh(\chi) E_1 Q_n=
Q_n E_{n}^* \Hh(\chi) E_n Q_n= E_n^* \Qh_1 \Hh(\chi) \Qh_1 E_n.
$$
Here we have used the invariance of $\Hh(\chi)$ under contraction. Hence we arrive at the sequence
\bqn\label{f.seq1}
E_n^*(Q_1+\beta Q_1\Hh(\chi) Q_1) E_n.
\eqn
Because $E_n$ is an isometry, it follows from Prop.~\ref{p22.1} that $E_n^* \Qh_1 \Hh(\chi) \Qh_1 E_n$
is self-adjoint and has spectrum contained in $[0,1]$. Thus, if $\beta\notin(-\iy,-1]$, then the sequence is stable.

Conversely, if the sequence (\ref{f.seq1}) is stable, then so is 
$$
E_n E_n^*(Q_1+\beta Q_1 H_\R(\chi) Q_1)E_n E_n^*.
$$
Now observe that $E_nE_n^*$ converges strongly to the identity operator on $L^2(\R_+)$ as $n\to \iy$.
It follows that $Q_1+\beta Q_1 H_\R(\chi) Q_1$ is invertible, which by Prop.~\ref{p22.1} implies
that $\beta\notin(-\iy,-1]$.
\end{proof}

\begin{theorem}\label{t3.7}
Let $a\in PC$ be continuous on $\T\setminus\{1,-1\}$.
Then the sequence $\{Q_n+Q_nH(a)Q_n\}_{n\in\Z_+}$ is stable if
$$\frac{a(1+0)-a(1-0)}{2i} \notin (-\iy,-1]
\quad\mbox{and}\quad
\frac{a(-1-0)-a(-1+0)}{2i} \notin (-\iy,-1].
$$ 
\end{theorem}
\begin{proof}
Let us write $a=a_1+a_2$, where $a_1$ is continuous on $\T\setminus\{1\}$ and 
$a_2$ is continuous on $T\setminus\{-1\}$. The previous proposition implies that
the sequences $Q_n+Q_nH(a_1)Q_n$ and $Q_n +Q_n H(a_3) Q_n$ are stable where
$a_3(t)=- a_2(-t)$. By applying the flip $f(t)\mapsto f(-t)$ (which is an isomorphism) to the left and the right hand side of the last sequence it follows that $Q_n +Q_n H(a_2) Q_n$ is stable.

Our next claim is that $Q_n H(a_1) Q_n H(a_2) Q_n$ tends to zero in the operator norm as $n\to\iy$.
To see this, let $V_{\pm n}=T(t^{\pm n})$. From (\ref{f.THabc}) we have
$Q_n=V_{n}V_{-n}$ and $V_{-n}H(b)=H(b t^{-n})= H(b) V_{n}$. Hence we can write
$$
Q_n H(a_1) Q_n H(a_2) Q_n=Q_n V_{-n} H(a_1)H(a_2) V_n Q_n.
$$
Notice that $V_{-n}\to 0 $ strongly and that the product $H(a_1)H(a_2)$ is a compact operator. The latter holds because we can write
\bqn
H(a_1)H(a_2) &=& H(a_1)T(f_1)H(a_2)+H(a_1)T(f_2)H(a_2)
\nn\\
&=& \Big( H(a_1f_1)-T(a_1)H(f_1)\Big)H(a_2) + H(a_1)\Big(H(f_2a_2)-H(f_2)T(\ta_2)\Big)
\nn
\eqn
with continuous even functions $f_1, f_2$ satisfying $f_1+f_2=1$ and such that $f_k$ is identically
zero on a neighborhood of the discontinuity of $a_k$. A similar argument shows that 
$Q_n H(a_2) Q_n H(a_1) Q_n$ converges to zero in the operator norm.

Now abbreviate $A_n=Q_nH(a_1)Q_n$ and $B_n=Q_n H(a_2) Q_n$. Then 
$$
Q_n+Q_nH(a)Q_n= Q_n+ A_n+B_n.
$$
We claim that $Q_n -(Q_n+A_n)\iv A_n -(Q_n+B_n)\iv B_n$ is an asymptotic inverse. If we multiply these two expressions we obtain
$$
Q_n -(Q_n+A_n)\iv A_n -(Q_n+B_n)\iv B_n+A_n -(Q_n+A_n)\iv A_n^2+B_n-(Q_n+B_n)\iv B_n^2
$$
plus terms tending to zero in the norm (because they contain the products $A_nB_n$ or $B_n A_n$). The last expression simplifies to $Q_n$.
This completes the proof.
\end{proof}

One can show that the previous result is actually an ``if and only if'' statement. Since we will not need
it we omit a proof. Our desired result of this section is the following.

\begin{corollary}\label{c3.8}
Let $\alpha,\beta\in\C$ be such  that $|\Re \alpha|<1$, $|\Re\beta|<1$, and let $c_+\in G\W_+$. Then the sequence
$$
\left\{P_n(I+H(\psi))\iv P_n\right\}_{n\in\Z_+} \quad \mbox{with}\quad  \psi=\tc_+ c_+\iv u_{-1/2-\alpha,1}u_{1/2+\beta,-1}
$$
is stable.
\end{corollary}
\begin{proof}
The invertibility of $I+H(\psi)$ on $H^2(\T)$ follows from Theorem \ref{t21}. Using Lemma \ref{l21.1} the stability can be reduced to the stability of $Q_n(I+H(\psi))Q_n$. We can apply the previous theorem, and for  this we compute 
$$
\frac{\psi(1+0)-\psi(1-0)}{2i}=\frac{e^{\pi i(1/2+\alpha)}-e^{-\pi i(1/2+\alpha)}}{2i}=\cos(\pi \alpha),
$$
$$
\frac{\psi(-1-0)-\psi(-1+0)}{2i}=\frac{e^{\pi i(1/2+\beta)}-e^{-\pi i(1/2+\beta)}}{2i}=\cos(\pi \beta).
$$
The corresponding conditions are fulfilled if and only if $\Re\alpha\notin 1+ 2\Z$ and $\Re\beta\notin 1+2\Z$.
\end{proof}

%%%%%%%%

\subsection{Stability results for operators with approximating symbols}
\label{s32}

In Sections \ref{s4} and \ref{s7} on we are going to approximate operators of the form
$$
(I+H(\phi))\iv
$$
with a certain symbol $\phi\in PC(\T)$ by operators of the same type
but with smooth symbols $\phi_\mu$. Clearly, we cannot expect approximation in the norm. What is sufficient for our purposes is the approximation in the strong operator topology. In view of Lemma \ref{l.basic}(iii) we need to examine the stability of operators $I+H(\phi_\mu)$, where $\mu\in[0,1)$ is the approximation parameter. The stability result is non-trivial, and we make in fact use of results established by the author and Silbermann in \cite{ES1}.

We will use the results of this subsection at two different places, namely in the proof of 
Theorem \ref{t3.2} and Proposition \ref{p7.4}.

Let us introduce the necessary notation. For $\mu\in[0,1)$ and $\tau\in\T$, define the composition operators $G_{\mu}$ and $Y_\tau$ acting on $L^\iy(\T)$ by
\be\label{f.Gm}
(G_{\mu} f)(t) = f\left(\frac{t+\mu}{1+\mu t}\right),\qquad (Y_\tau f)(t)=f(\tau t),
\ee
and the operators $R_{\mu}$ acting on $L^2(\T)$ by
\be\label{f.Rm}
(R_{\mu} f)(t) = \frac{\sqrt{1-\mu^2}}{1+\mu t} f\left(\frac{t+\mu}{1+\mu t}\right).
\ee
The operator $Y_\tau$ also acts on $L^2(\T)$, and both $R_{\mu}$ and $Y_\tau$ are unitary operators. The Hardy space $H^2(\T)$ is an invariant subspace of $R_{\mu}$ and $Y_\tau$ as well as of their adjoints. We will use the same notation for the restriction of $R_{\mu}$ and $Y_\tau$ onto $H^2(\T)$. These restrictions are unitary operator on $H^2(\T)$, too. Moreover, for  $\phi\in L^\iy(\T)$ we have
\bqn\label{RGH}
R_{\mu} H(\phi) R_{\mu}^* =  H(G_{\mu}\phi),\qquad
R_{\mu} T(\phi) R_{\mu}^* =  T(G_{\mu}\phi),\qquad \mu\in[0,1),
\eqn
and
\bqn\label{RGH2}
Y_\tau H(\phi) Y_\tau^* = \tau H(Y_\tau\phi)\quad \mbox{if }\tau=\pm1,\qquad
Y_\tau T(\phi) Y_\tau^* = T(Y_\tau \phi) \quad \mbox{if }\tau\in\T .
\eqn
The just mentioned statements are easy to prove (see also \cite[Sec.~5.1]{ES1}).

\begin{theorem}\label{t.23st}
Let $\alpha,\beta\in \C$ be such that $|\Re\alpha|<1$ and $|\Re\beta|<1$,  let 
$c_+\in G\W_+$, and let
$$
a_\mu(t)=\left(\frac{1-\mu t}{1-\mu t\iv}\right)^{-1/2-\alpha},\qquad 
b_\mu(t)=\left(\frac{1+\mu t}{1+\mu t\iv}\right)^{1/2+\beta},\qquad \mu\in[0,1).
$$
Then the (generalized) sequence of operators
$$
\Big\{I+H(\psi_\mu)\Big\}_{\mu\in[0,1)}, \qquad \psi_\mu:= \tilde{c}_+ c_+\iv a_\mu b_\mu
$$
is stable on $H^2(\T)$.
\end{theorem}
\begin{proof}
In order to prove the stability we apply the results of
\cite[Secs.~4.1--4.2]{ES1}.  These results establish the existence of certain
mappings $\Phi_0$ and $\Phi_\tau$, $\tau\in\T$, which are defined as
$$
\Phi_0[\psi_\mu ]:= \mu\mbox{-}\lim_{\mu \to1}  \psi_\mu ,\qquad
\Phi_\tau[\psi_\mu ]:= \mu\mbox{-}\lim_{\mu \to1} G_\mu  Y_{\tau}\psi_\mu ,
$$
where $\mu\mbox{-}\lim$ stands for the limit in measure. 
It is easy to see that 
$$
\Phi_0[a_\mu ] =u_{-1/2-\alpha,1},\qquad
\Phi_0[b_\mu ] =u_{1/2+\beta,-1},
$$
and 
$$
\Phi_\tau[a_\mu ]=  u_{-1/2-\alpha,1}(\tau),\quad \mbox{ if }\tau\neq1,\qquad
\Phi_{1}[a_\mu ]=u_{1/2+\alpha,-1},
$$
and
$$
\Phi_\tau[b_\mu ] =u_{1/2+\beta,-1}(\tau),\quad \mbox{ if }\tau\neq-1,\qquad
\Phi_{-1}[b_\mu ]=u_{-1/2-\beta,-1}.
$$
Since $\psi_\mu  =\tc_+ c_+\iv a_\mu  b_\mu $ we conclude
\bqn
\Phi_0[\psi_\mu ] &=& \tc_+ c_+\iv u_{-1/2-\alpha,1}u_{1/2+\beta,-1},\nn\\[1ex]
\Phi_1[\psi_\mu ] &=& u_{1/2+\alpha,-1},\nn\\[1ex]
\Phi_{-1}[\psi_\mu ] &=& u_{-1/2-\beta,-1},\nn\\[1ex]
\Phi_\tau[\psi_\mu ] &=& \mbox{constant function, }\quad\tau\in \T\setminus\{-1,1\}.\nn
\eqn
The stability criterion in \cite{ES1} (Thm.~4.2 and Thm.~4.3) says that $\{I+H(\psi_\mu)\}_{\mu\in[0,1)}$ 
is stable if and only if the operators
\begin{itemize}
\item[(i)] $\Psi_0[I+H(\psi_\mu )]=I+H(\Phi_0[\psi_\mu ])=I+H( \tc_+ c_+\iv u_{-1/2-\alpha,1}u_{1/2+\beta,-1})$,
\item[(ii)] $\Psi_1[I+H(\psi_\mu )]=I+H(\Phi_1[\psi_\mu ])=I+H(u_{1/2+\alpha,-1})$,
\item[(iii)] $\Psi_{-1}[I+H(\psi_\mu )]=I-H(\Phi_{-1}[\psi_\mu ])=I-H(u_{-1/2-\beta,-1})$,
\item[(iv)] $\Psi_\tau[I+H(\psi_\mu )]=$
$$
\twomat{I&0\\0&I}+\twomat{P&0\\0&Q}\twomat{M(\Phi_\tau[\psi_\mu ])&0\\0&
M(\widetilde{\Phi_{\bar{\tau}}[\psi_\mu ]})}\twomat{0&I\\ I&0}
\twomat{P&0\\0&Q}=\twomat{I&0\\0&I}
$$
($\tau\in\T$, $\mathrm{Im}(\tau)>0$)
\end{itemize}
are invertible. 

The operator $I-H(u_{-1/2-\beta,-1})$ is invertible if and only if
$I+H(u_{-1/2-\beta,1})$ is invertible because these two operators can be related to each other 
by a rotation operator $Y_{-1}:f(t)\mapsto f(-t)$ acting on $H^2(\T)$ (see (\ref{RGH2})).
Now the invertibility of the operators (i)--(iii) follows from Theorem \ref{t21}.
\end{proof}

%%%%%%%%%%%%

%%%%%%%%%%%%%%%%%%%%%%%%%%%%%%%%%%%%

\section{A formula for Hankel determinants}
\label{s4}

For $b\in L^1[-1,1]$, let $H_n[b]$ stand for the $n\times n$ Hankel matrix
\be\label{f.Hdet}
H_n[b]=(b_{j+k+1})_{j,k=0}^{n-1},
\ee
where
\be\label{f.moments}
b_{k}=\frac{1}{\pi}\int_{-1}^1 b(x) (2x)^{k-1}\, dx
\ee
are the (scaled) moments of $b$.
In this section we are going to derive a formula of the type
$$
\det H_n[b]= G^n\det\Big( P_n(I+H(\psi))\iv P_n \Big),
$$
where $b\in L^1[-1,1]$ is function of the form 
$b(x)=(1-x)^\alpha(1+x)^\beta b_0(x)$ and $b_0(x)$ is a sufficiently smooth and 
nonvanishing function on $[-1,1]$. The constant $G$ and the function $\psi\in PC(\T)$ depend on
the function $b$.
This formula will allow us in the next section (Sec.~\ref{s5}) to express the determinant of a Hankel matrix as a determinant of the type appearing on the right hand side. The invertibility of $I+H(\psi)$
will be guaranteed by Theorem \ref{t21}. The $P_n$'s are the finite sections (\ref{f.Pn}).
The above formula is a generalization of formulas of the same type established in \cite{E1, E2} for particular values of $\alpha, \beta$.

Recall the definition of the Wiener algebra $\W$ given in Sec.~\ref{s22}. A function $a\in\W$ is said to admit a canonical Wiener-Hopf factorization in $\W$ 
if it can be represented in the form
\be
a(t)=a_-(t)a_+(t),\qquad t\in\T,
\ee
where $a_\pm\in G\W_\pm$. It is well known (see, e.g., \cite{BS}) that $a\in\W$ 
admits a canonical Wiener-Hopf factorization in $\W$ if and only if $a\in G\W$ and if 
the winding number of $a$ is zero. This, in turn, is equivalent to the condition that 
$a$ possesses a logarithm $\log a\in \W$. In this case, one can define the
geometric mean
\bqn
G[a] &:=& \exp\Big(\frac{1}{2\pi}\int_{0}^{2\pi} \log a(e^{i\theta})\, d\theta\Big).
\eqn
This definition does not depend on the particular choice of the logarithm.

We will assume that $a\in G\W$ is an even function. Then $a$ has winding number zero
and thus possesses a canonical Wiener-Hopf factorization. Moreover, the factors are related to each other by
$a_-=\gamma \ta_+$ with some nonzero constant $\gamma$.

The following theorem is cited from \cite[Thm.~4.5]{E1}. It is the immediate consequence of two other results, namely, Thm.~2.3 of \cite{BE2} and Prop.~3.9 of \cite{BE3}. To give some, but not all details, we remark that the last two mentioned results establish the identites
$$
\det H_n[b]=\det \Big( P_n(T(a)+H(a))P_n\Big)
$$
and 
$$
\det \Big( P_n(T(a)+H(a))P_n\Big)=G[a]^n\det\Big( P_n (I+H(\psi))\iv P_n\Big).
$$
Therein it is necessary to assume that the symbols are smooth (or, more precisely, that the Wiener-Hopf factors of $a$ are bounded). We remark that the invertibility of $I+H(\psi)$ is guaranteed again by
Theorem \ref{t21}, which can be applied with the jump functions being absent.
A more elementary argument for the invertibility of $I+H(\psi)$ with (certain) smooth symbols
$\psi$ and an explicit formula for the inverse is given in Prop.~3.9 of \cite{BE3}.

\begin{theorem}\label{t3.1}
Let $a\in G\W$ be an even function which possesses a Wiener-Hopf factorization
$a(t)=a_-(t)a_+(t)$. Define $\psi(t)=\ta_+(t)a_+\iv(t)$ and
\bqn\label{f.TH-H} 
b(\cos \theta) = a(e^{i\theta}) \sqrt{\frac{1+\cos \theta}{1-\cos\theta}}.
\eqn
Then $I+H(\psi)$ is invertible on $H^2(\T)$ and 
\bqn
\det H_n[b] &=&G[a]^n \det\Big( P_n (I+H(\psi))\iv P_n\Big).
\eqn
\end{theorem}

The next theorem, which is the main result of this section, follows formally from the previous one with 
$$
a(e^{i\theta})=c(e^{i\theta})(2-2\cos\theta)^{1/2+\alpha}(2+\cos\theta)^{-1/2+\beta},
$$
$$
a_+(t)=c_+(t) (1-t)^{1/2+\alpha}(1+t)^{-1/2+\beta}, \quad 
a_-(t)=c_-(t)(1-t\iv)^{1/2+\alpha}(1+t\iv)^{-1/2+\beta}
$$
if one considers the generalized Wiener-Hopf factorizations
$$
u_{-1/2-\alpha,1}(t)=(1-t)^{-1/2-\alpha}(1-t\iv)^{1/2+\alpha},\quad
u_{1/2-\beta,-1}(t)=(1+t)^{1/2-\beta}(1+t\iv)^{-1/2+\beta}.
$$
Of course, therein the Wiener-Hopf factors and the inverses are not bounded.
In fact, in order to make the argument precise, we have to use an approximation argument,
and we have to apply the stability results of Section \ref{s32}.
 
\begin{theorem}\label{t3.2}
Let $c\in G\W$ be an even function which possesses a Wiener-Hopf factorization
$c(t)=c_-(t)c_+(t)$. Define 
\bqn
\psi(t) &=& \tc_+(t)c_+\iv(t) u_{-1/2-\alpha,1}(t)u_{1/2-\beta,-1}(t) ,\qquad t\in \T,
\\[1ex]
b(\cos \theta) & =& c(e^{i\theta})(2-2\cos \theta)^\alpha (2+2\cos \theta)^\beta,
\eqn
and assume that $|\Re\alpha|<1$, $|\Re\beta|<1$. Then 
$$
\det H_n[b]= G[c]^n\det\Big(P_n(I+H(\psi))\iv P_n\Big).
$$
\end{theorem}
\begin{proof}
For $\mu\in[0,1)$ we define the even functions
$$
a_\mu(t)= c(t)
\Big((1-\mu t)(1-\mu t\iv)\Big)^{1/2+\alpha}
\Big((1+\mu t)(1+\mu t\iv)\Big)^{-1/2+\beta},\qquad t\in\T.
$$
By formula (\ref{f.TH-H}), these functions $a_\mu $ give rise to functions
$$
b_\mu (\cos\theta)=c(e^{i\theta})\cdot
\frac{(1+\mu ^2-2\mu \cos\theta)^{1/2+\alpha}}{(2-2\cos\theta)^{1/2}}\cdot
\frac{(1+\mu ^2+2\mu \cos\theta)^{-1/2+\beta}}{(2+2\cos\theta)^{-1/2}}.
$$
Because $b_\mu \to b$ in the norm of $L^1[-1,1]$ as $\mu \to1$, 
it follows that (for fixed $n$)
$$
H_n[b]=\lim_{\mu \to1} H_n[b_\mu ].
$$
Furthermore, we observe that the canonical Wiener-Hopf factorization of $a_\mu $ is given
by $a_\mu (t)=a_{\mu ,-}(t) a_{\mu ,+}(t)$ with
\bqn
a_{\mu ,-}(t) &=& c_-(t)(1-\mu t\iv)^{1/2+\alpha}(1+\mu t\iv )^{-1/2+\beta},\nn
\\
a_{\mu ,+}(t) &=&c_+(t)(1-\mu t)^{1/2+\alpha}(1+\mu t)^{-1/2+\beta}.\nn
\eqn
Then 
$$
\psi_\mu (t)=\ta_{\mu ,+}(t)a_{\mu ,+}\iv(t)=
\tc_{+}(t)c_{+}\iv(t)
\left(\frac{1-\mu t}{1-\mu t\iv}\right)^{-1/2-\alpha}
\left(\frac{1+\mu t}{1+\mu t\iv}\right)^{1/2-\beta}.
$$
Using the fact that $G[a]=G[c]$ and employing Theorem \ref{t3.1} it follows that 
$$
\det H_n[b_\mu ]=G[c]^n\det\Big ( P_n ( I+H(\psi_\mu ))\iv P_n\Big ),
$$
whence
$$
\det H_n[b]= G[c]^n\lim_{\mu \to 1} \det\Big (P_n ( I+H(\psi_\mu ))\iv P_n\Big).
$$

Because
$$
\left(\frac{1-\mu t}{1-\mu t\iv}\right)^{-1/2-\alpha}
\to u_{- 1/2-\alpha,1}(t),\qquad
 \left(\frac{1+\mu t}{1+\mu t\iv}\right)^{1/2-\beta}
\to u_{1/2-\beta,-1}(t),
\qquad
$$
in measure as $\mu \to 1$,
it follows that $\psi_\mu \to \psi$ in measure as $\mu \to1$.
Hence, by Lemma \ref{l3.3},  $H(\psi_\mu)\to H(\psi)$ strongly on $H^2(\T)$ as $\mu \to1$.

Now Theorem \ref{t.23st} implies that the sequence $\{I+H(\psi_\mu )\}_{\mu\in[0,1)}$ is stable, while
Theorem \ref{t21} implies that the operator $I+H(\psi)$ is invertible on $H^2(\T)$.
Hence (see Lemma \ref{l.basic}(iii))
$$
(I+H(\psi_\mu ))\iv \to (I+H(\psi))\iv
$$
strongly on $H^2(\T)$ as $\mu \to1$. This allows us to conclude that 
(for fixed $n$) the matrices $P_n(I+H(\psi_\mu ))\iv P_n$ converge to 
$P_n (I+H(\psi))\iv P_n$ as $\mu \to 1$.  Consequently, their determinants also converge.
This completes the proof.
\end{proof}

We conjecture that the conditions on the parameters $\alpha$ and $\beta$ in the previous
theorem cannot be weakened. This is despite of the fact that the Hankel determinant
$H_n[b]$ is well defined if $\Re\alpha>-1$ and $\Re\beta>-1$, and  that the 
inverse of $I+H(\psi)$ exists under the same condition and the extra condition that
$\Re\alpha\notin 1+2\Z$ and $\Re\beta\notin 1+2\Z$. (For the last statement see 
\cite[Sec.~3.2]{BE3} or  \cite{BEp}.)

%%%%%%%%%%%%%%%%%%%%%%%%%%%%%%%%%%%%%%%%%%%%%%%%%%%%%%%%%%%%%%%%%%%%%%%

\section{The approximation of the Bessel kernel}
\label{s5}

After having established the necessary auxiliary results in the previous sections we are able to start 
with analyzing the quantity 
$$
P^{(\alpha)}(R) = \det (P_R-P_RB_\alpha P_R)|_{L^2[0,R]}.
$$
As outlined in the introduction our first goal (Step 1) is to establish an alternative expression for this quantity, and this will be accomplished at the end of the next section.

In this section we are going to approximate the Bessel kernel by the Jacobi kernel. 
This approximation corresponds to the fact that correlation functions for the Jacobi Unitary Ensemble yield in the hard edge scaling limit the Bessel kernel.

It is possible to consider the Bessel kernel $B_\alpha(x,y)$ for complex parameters $\alpha$ with $\Re\alpha>-1$, and the Jacobi kernel for complex parameters $\alpha,\beta$ with $\Re\alpha>-1$, $\Re\beta>-1$. Recall that the Jacobi kernel is defined by
\be\label{f.44}
K_n^{(\alpha,\beta)}(x,y) = \sqrt{w(x)w(y)}\; \sum_{k=0}^{n-1} p_k(x) p_k(y),\qquad -1< x,y < 1,
\ee
where $w(x)=(1-x)^\alpha(1+x)^\beta$ is the underlying weight, and $p_k(x)$ are the normalized orthogonal polynomials, i.e., 
\be\label{f.44b}
\int_{-1}^1p_j(x)p_k(x)w(x)\, dx =\delta_{jk},
\ee
which are, of course, up to a constant equal to the Jacobi polynomials.
We define the rescaled kernel as 
$$
\hat{K}^{(\alpha,\beta)}_n(x,y)=\frac{\sqrt{xy}}{n^2} K_n^{(\alpha,\beta)}(1-\frac{x^2}{2n^2},1-\frac{y^2}{2n^2}), \qquad 0<x,y < 2 n.
$$
For sufficiently large $n$ (namely, $2n>R$) we can think of $\hat{K}_{n}^{(\alpha,\beta)}$ 
as an integral operator acting on $L^2[0,R]$, where $R>0$ is considered fixed.

From (\ref{f.44}) it is clear that the operators $K_n^{(\alpha,\beta)}$ and
$\hat{K}_n^{(\alpha,\beta)}$ are finite rank operators with rank equal to $n$ (or at most $n$).

\begin{theorem}\label{t4.1}
Let $R>0$ be fixed.
\begin{itemize}
\item[(i)]
For $\Re \alpha>-1$ the operator $P_R B_\alpha P_R|_{L^2[0,R]}$ is trace class, and hence the operator determinant
$$
P^{(\alpha)}(R)= \det (P_R-P_R B_\alpha P_R)|_{L^2[0,R]}
$$
is well defined. Moreover, the mapping 
$\alpha\mapsto P_R B_\alpha P_R|_{L^2[0,R]}$ is an analytic (trace class) operator-valued function.
\item[(ii)]
For $\alpha, \beta>-1$ being real, we have $P_R \hat{K}_n^{(\alpha,\beta)} P_R\to P_R B_\alpha P_R$
in the trace norm and hence
\be\label{P.lim}
P^{(\alpha)}(R)=\lim_{n\to\iy} \det(P_R-P_R\hat{K}_n^{(\alpha,\beta)}P_R).
\ee
\end{itemize}
\end{theorem}
\begin{proof}
(i): The Bessel operator $B_\alpha$ is by definition the integral operator with the kernel 
(\ref{Bessel.kernel}). Because this kernel can also be expressed by (\ref{Bessel.int2}),
we can write $P_RB_\alpha P_R=X_{R,\alpha} Y_{R,\alpha}$ where 
$$
X_{R,\alpha}:L^2[0,1]\to L^2[0,R],\qquad  Y_{R,\alpha}:L^2[0,R]\to L^2[0,1]
$$ 
are integral operators with the kernels
$$
X_{R,\alpha}(x,t)=\sqrt{tx}J_\alpha(tx),\qquad Y_{R,\alpha}(t,y)=\sqrt{ty} J_\alpha(ty).
$$
Using the asymptotics of the Bessel function at zero, $J_\alpha(t)\sim (t/2)^\alpha/ \Gamma(1+\alpha)$,
it is easy to conclude that both $X_{R,\alpha}$ and $Y_{R,\alpha}$ are Hilbert-Schmidt operators (see \cite[Thm.~2.11]{S}).
Hence $P_RB_\alpha P_R=X_{R,\alpha} Y_{R,\alpha}$ is trace class.

It is a basic fact that the function $z^{-\alpha} J_\alpha(z)=:g(\alpha,z)$ is entire in both $\alpha$ and $z$
(see \cite[Sec.\ II,\ 9.3]{O}). Thus, differentiating $\sqrt{tx}J_\alpha(tx)=(tx)^{\alpha+1/2}g(\alpha,tx)$ with respect to $\alpha$ yields $\ln(tx) (tx)^{\alpha+1/2}g(\alpha,tx)+(tx)^{\alpha+1/2}(\partial_\alpha g)(\alpha,tx)$.
This implies that the derivative (with respect to complex $\alpha$) of the operator-valued function $\alpha\mapsto X_{R,\alpha}$  exists and is in fact a Hilbert-Schmidt. The same holds for the derivative of $Y_{R,\alpha}$. Hence $P_R B_\alpha P_R$ has a complex derivative with respect to $\alpha$, which is trace class.

(ii): This part was proved by Kuijlaars and Vanlessen \cite{KV02}. They actually established the result 
not just for the usual Jacobi weight, but for a modified Jacobi weight. Adapting their notation to ours
they show, using the Riemann-Hilbert method, that 
\be\label{f.est1}
\hat{K}^{(\alpha,\beta)}_n(x,y)= B_\alpha(x,y)+O\left(\frac{x^{\alpha+1/2}y^{\alpha+1/2}}{n}\right)\quad 
\mbox{ as } n\to\iy,
\ee
where the error term holds uniformly in $x$ and $y$ on bounded subsets of $(0,\iy)$.
From the decomposition $P_RB_\alpha P_R=X_{R,\alpha} Y_{R,\alpha}$ it also follows easily that
$|B_\alpha(x,y)|\le C x^{\alpha+1/2} y^{\alpha+1/2}$ uniformly in $x$ and $y$ on bounded subsets of $(0,\iy)$. These estimates imply that $P_R\hat{K}_n^{(\alpha,\beta)}P_R\to P_RB_\alpha P_R$ weakly
and $\tr(P_R\hat{K}_n^{(\alpha,\beta)}P_R)\to\tr( P_RB_\alpha P_R)$ as $n\to\iy$. Because
for real $\alpha$ and $\beta$, the operators $\hat{K}_n^{(\alpha,\beta)}$ and $B_\alpha$ are positive,
using \cite[Thm.~2.20]{S}, one can conclude that $P_R\hat{K}_n^{(\alpha,\beta)}P_R\to P_RB_\alpha P_R$ 
in the trace norm.
\end{proof}

The question arises whether the restriction to real $\alpha,\beta$ in (ii) is necessary. We conjecture that it is not. It seems likely that one can prove estimate (\ref{f.est1}) by repeating the arguments of \cite{KV02}. Moreover, one could perhaps also prove (\ref{f.est1}) by using the classical results on the 
asymptotics of the Jacobi polynomials.

Assuming the validity of this estimate it follows immediately that the traces converge and 
that $P_R\hat{K}_n^{(\alpha,\beta)}P_R\to P_RB_\alpha P_R$ in the Hilbert-Schmidt norm. 
Using regularized determinants the validity of (\ref{P.lim}) would follow for complex $\alpha$.

It is probably true that even $P_R\hat{K}_n^{(\alpha,\beta)}P_R\to P_RB_\alpha P_R$ in the
trace norm for complex $\alpha$. Since the operators are not positive definite 
the argument would be more complicated. This statement is dispensable for
our purposes, and hence we will not discuss this issue further.

The restriction to real $\alpha,\beta$ is only temporary. It will be removed at the end of Section 
\ref{s6} using an analyticity argument.

\begin{proposition}\label{p5.2}
Let $R>0$ be fixed, and assume that $\alpha,\beta>-1$ are real. Then
\be\label{f.47}
P^{(\alpha)}(R)=
\exp\Big(-\frac{R^2}{4}+\alpha R\Big)
\lim_{n\to\iy}
\frac{\det\left(\int_{-1}^1x^{j+k} \hat{w}_n(x)\, dx\right)_{j,k=0}^{n-1}}
{\det\left(\int_{-1}^1x^{j+k} \hat{w}(x)\, dx\right)_{j,k=0}^{n-1}}
\ee
with
$$
\rho_n=1-\frac{R^2}{4 n^2},\quad \mu_n=\frac{2-\rho_n-2\sqrt{1-\rho_n}}{\rho_n}=1-\frac{R}{n}+O(n^{-2}).
$$
and
\be\label{weights}
\hat{w}_n(x)=(1+\mu_n^2-2\mu_n x)^\alpha (2+2x)^\beta,\qquad
\hat{w}(x)=(2-2x)^\alpha(2+2x)^\beta.
\ee
\end{proposition}
\begin{proof}
Suppose that $2n>R$. From (\ref{f.44}) and (\ref{f.44b}) it follows that 
we can decompose $\hat{K}_n^{(\alpha,\beta)}|_{L^2[0,R]}=AB$, where
$$
A: \{x_k\}_{k=0}^{n-1}\in\C^n\mapsto f\in L^2[0,R], \quad f(x)=\frac{\sqrt{x}}{n} \sqrt{w(1-\frac{x^2}{2n^2})}\sum_{k=0}^{n-1}p_k(1-\frac{x^2}{2n^2}) x_k
$$
and
$$
B: f\in L^2[0,R]\mapsto \{y_k\}_{k=0}^{n-1}\in \C^n, \quad y_k=\int_{0}^R \frac{\sqrt{y}}{n}\sqrt{w(1-\frac{y^2}{2n^2})} p_k(1-\frac{y^2}{2n^2}) f(y)\, dy.
$$
Recall that $w(x)=(1-x)^\alpha (1+x)^\beta$.
Then $BA$ is an $n\times n$ matrix
\bqn
BA &=& 
\left(\int_{0}^R p_j(1-\frac{x^2}{2n^2})p_k(1-\frac{x^2}{2n^2}) w(1-\frac{x^2}{2n^2})\,\frac{x}{n^2} \,dx\right)_{j,k=0}^{n-1}
\nn\\
&=&
\left(\int_{1-\frac{R^2}{2n^2}}^1 p_{j}(y)p_k(y)w(y)\, dy\right)_{j,k=0}^{n-1}.
\nn
\eqn
Hence (using the orthonormality of the $p_k$'s and the formula $\det(I-AB)=\det(I-BA)$)
\bqn
\det(I-\hat{K}_n^{(\alpha,\beta)})|_{L^2[0,R]}
&=&
\det\left(\int_{-1}^{1-\frac{R^2}{2n^2}}p_{j}(y)p_k(y)w(y)\, dy\right)_{j,k=0}^{n-1}.\nn
\eqn
This equals 
\bqn%\label{f.int1}
\left(\prod_{k=0}^{n-1} \sigma_k^2\right) \det \left(\int_{-1}^{1-\frac{R^2}{2n^2}} x^{j+k}w(x)\, dx\right)_{j,k=0}^{n-1}\nn
\eqn
where $\sigma_k$ is the coefficient of the leading term of $p_k(x)$.
Obviously,
$$
1=\left(\prod_{k=0}^{n-1} \sigma_k^2\right) \det \left(\int_{-1}^{1} x^{j+k}w(x)\, dx\right)_{j,k=0}^{n-1}.
$$
{}From this we can eliminate the product of the $\sigma_k^2$'s, and we arrive at
\bqn\label{f.int1}
\det(I-\hat{K}_n^{(\alpha,\beta)})|_{L^2[0,R]}
=
\frac{\det\left(\int_{-1}^{1-\frac{R^2}{2n^2}} x^{j+k} w(x)\, dx\right)_{j,k=0}^{n-1}}
{\det\left(\int_{-1}^1x^{j+k} w(x)\, dx\right)_{j,k=0}^{n-1}}.
\eqn
In the integral appearing in the numerator we make a substitution $x=\rho_n(y+1)-1$ with $\rho_n=1-\frac{R^2}{4 n^2}$
to obtain 
$$
\rho_n^{j+k+1}\int_{-1}^1(y+\gamma_n)^{j+k}w(\rho_n(y+1)-1)\, dy
$$
$\gamma_n=1-1/\rho_n$. By expanding $(y+\gamma_n)^j(y+\gamma_n)^k$ using twice the binomial formula and by performing row and column operations, the corresponding determinant
equals
$$
\rho_n^{n^2} \det\left(\int_{-1}^{1} y^{j+k} w(\rho_n(y+1)-1))\, dy\right)_{j,k=0}^{n-1}.
$$
The new weight therein evaluates to
$$
(2-\rho_n-\rho_ny)^\alpha \rho_n^\beta (y+1)^\beta=
((1+\mu_n^2)/2-\mu_n y)^\alpha (y + 1)^\beta \rho_n^{\alpha+\beta}\mu_n^{-\alpha}  
$$
where $(1+\mu_n^2)/(2\mu_n)=(2-\rho_n)/\rho_n$ and $\mu_n\in (0,1)$.
Thus the determinant in the numerator of (\ref{f.int1}) equals
$$
\rho_n^{n^2}\rho_n^{n(\alpha+\beta)}\mu_n^{-n\alpha} \det \left(\int_{-1}^{1} x^{j+k}
((1+\mu_n^2)/2-\mu_n x)^\alpha (x + 1)^\beta\, dx\right)_{j,k=0}^{n-1}.
$$
This is almost the same expression as in the numerator of (\ref{f.47}) except for a factors of ``2''
in the weight. As the denominator in (\ref{f.int1}) with the weight $w(x)$ lacks the same factor
as compared to the weights $\hat{w}(x)$, the corresponding quotient is the same.
Now it remains to remark that the limit of $\rho_n^{n^2}\rho_n^{n(\alpha+\beta)}\mu_n^{-n\alpha} $
as $n\to\iy$ gives the constant $\exp(-R^2/4+\alpha R)$.
This completes the proof.
\end{proof}

We could, of course, compute the determinant appearing in the denominator of (\ref{f.47}) explicitly.
It is, up to a factor $2^{(\alpha+\beta)n}$, just the product of the $\sigma_k^{-2}$'s appearing in the proof, which are known quantities for the Jacobi polynomials. We could see that the product is nonzero,
and hence the determinant in the denominator of (\ref{f.47}) is nonzero, too.
However, this computation will not be of help for us since we need this expression in order to perform some kind of cancellation with the determinant in the numerator.

\begin{theorem}\label{t4.3}
Let $R>0$ be fixed, and assume that $-1<\alpha,\beta <1$ are real. Then
\bqn\label{f.16}
P^{(\alpha)}(R) &=&
\exp\Big(-\frac{R^2}{4}+\alpha R \Big)\lim_{n\to\iy}
\frac{\det \Big( P_n(I+H(\psi_n))\iv P_n\Big)}{\det\Big( P_n(I+H(\psi))\iv P_n\Big)}
\eqn
where 
\bqn
\psi_n(t) &=&   \left( \frac{1-\mu_nt}{1-\mu_n t\iv}\right)^{-\alpha} u_{-1/2,1}(t) u_{1/2-\beta,-1}(t),
\nn\\[1ex]
\psi(t) &=& u_{-1/2-\alpha,1}(t) u_{1/2-\beta,-1}(t),\qquad t\in\T,\nn
\eqn
and $\mu_n\in[0,1)$ is a sequence such that $\mu_n=1-R/n+O(1/n^2)$ as $n\to\iy$.
The expression under the limit in (\ref{f.16}) is well defined for each $n$.
\end{theorem}
\begin{proof}
Using the notation (\ref{f.Hdet}) for the Hankel determinant and the weights (\ref{weights})
it is easy to see that 
$$
\frac{\det\left(\int_{-1}^1x^{j+k} \hat{w}_n(x)\, dx\right)_{j,k=0}^{n-1}}
{\det\left(\int_{-1}^1x^{j+k} \hat{w}(x)\, dx\right)_{j,k=0}^{n-1}}
=\frac{\det H_n[\hat{w}_n]}{\det H_n[\hat{w}]}.
$$
Notice that a common factor due to the particular scaling of the moments (\ref{f.moments}) cancels out.
Now we only need to apply Theorem \ref{t3.2} twice in order to express the two Hankel determinants
as determinants of the other type.

We also remark that the operators $I+H(\psi_n)$ and $I+H(\psi)$ are invertible due to Theorem \ref{t21}.
Moreover, the determinant $\det\Big( P_n(I+H(\psi))\iv P_n\Big)$ is nonzero for all $n$.
This can be seen by relating it to the determinant in the denominator of (\ref{f.47}), which  is nonzero as
has been pointed out above.
\end{proof}

%%%%%%%%%%%%%%%%%%%%%%%%%%%%%%%%%%%%%%%%%%%%%%%%%%%%%%%%%%%%%%%%%%%%%%%

\section{Some asymptotic analysis}
\label{s6}

The goal now is to identify the limit of
$$
\frac{\det\Big( P_n(I+H(\psi_n))\iv P_n\Big)}{\det\Big(  P_n(I+H(\psi))\iv P_n\Big)}
$$
as $n\to\iy$. Although the functions $\psi_n$ and $\psi$ depend on $\beta$, in view of (\ref{f.16})
the limit is independent. Hence  it is sufficient to do the analysis for one particular value of $\beta$. It turns out that the choice $\beta=-\alpha$ makes things sufficiently simple. It is also possible to do the
analysis for arbitrary $\beta$, but the resulting formula is the same.

In case $\beta=-\alpha$, the functions $\psi_n$ and $\psi$ defined in Theorem \ref{t4.3}
become
\be\label{f.17}
\psi(e^{i\theta})= u_{-1/2-\alpha,1}(e^{i\theta}) u_{1/2+\alpha,-1}(e^{i\theta})=
\left\{\ba{ll} e^{i(\alpha+1/2)\pi} & \mbox{ if } 0<\theta<\pi\\
e^{-i(\alpha+1/2)\pi} & \mbox{ if } -\pi<\theta<0
\ea\right.
\ee
and
\be\label{f.17b}
\psi_n(t) =   \left( \frac{1-\mu_nt}{1-\mu_n t\iv}\right)^{-\alpha} u_{-1/2,1}(t) u_{1/2+\alpha,-1}(t).
\ee
For the purpose of the following lemma, recall the definition of the operators
$G_\mu$ and $R_\mu$ given in (\ref{f.Gm}) and (\ref{f.Rm}). Clearly, $\psi$ is invariant under 
$G_\mu$, i.e., $G_\mu \psi=\psi$. In order to rewrite $\psi_n$ notice that 
$$
\left(\frac{1-\mu_n t}{1-\mu_n t\iv}\right)^{-\alpha}=
 (-t)^{-\alpha} \left(-\frac{t-\mu_n}{1-\mu_n t}\right)^\alpha=
u_{-\alpha,1}(t)\cdot (G_{\mu_n}\iv u_{\alpha,1})(t).
$$
Here the principle values of the power functions are considered.
From this we conclude 
$\psi_n = \psi\cdot  (G_{\mu_n}\iv u_{\alpha,1})$. 
Now we define
\be\label{psihat}
\hat{\psi}:=G_{\mu_n }\psi_n = \psi \cdot u_{\alpha,1}=u_{-1/2,1} \cdot u_{1/2+\alpha,-1}.
\ee
Remark that the sequence $\mu_n$ used in these above formulas satisfies
\bqn\label{mu.n}
\mu_n = 1-\frac{R}{n} + O(n^{-2}),\qquad n\to\iy,
\eqn
where $R>0$ is fixed. Finally, let us define the function
$$
h_{R}(t)=\exp\Big (R\,\frac{t-1}{t+1}\Big), \qquad t\in\T,
$$
and the operator $\Pi_R=H(h_R)^2$ acting on $H^2(\T)$.

\begin{lemma}\label{l6.1}
Let $R>0$ be fixed, $|\Re\alpha|<1$, and let $\psi,\psi_n,$ and $\hat{\psi}$ be defined by
(\ref{f.17}), (\ref{f.17b}), (\ref{psihat}), and (\ref{mu.n}). Then the following holds.
\begin{enumerate}
\item[(a)] 
The operators $I+H(\psi)$ and $I+H(\hat{\psi})$ are invertible on $H^2(\T)$, and the sequence
$\{I+H(\psi_n)\}_{n\in\Z_+}$ is stable.
\item[(b)]
We have, as $n\to\iy$,
\be\label{f.19a}
R_{\mu_n}\Big(P_n(I+H(\psi ))\iv P_n\Big) R_{\mu_n}^* +R_{\mu_n}Q_n R_{\mu_n}^*\to \Pi_R(I+H(\hat{\psi})\iv )\Pi_R+(I-\Pi_R)
\ee
strongly and 
\be\label{f.19b}
R_{\mu_n} P_n\left( (I+H(\psi_n))\iv -(I+H(\psi))\iv \right) P_nR_{\mu_n}^*
\to
\Pi_R\left( (I+H(\hat{\psi}))\iv-(I+H(\psi))\iv\right) \Pi_R
\ee
in trace norm.
\end{enumerate}
\end{lemma}
\begin{proof}
(a):\  The invertibility of $I+H(\psi)$ and $I+H(\hat{\psi})$ follows from 
Theorem \ref{t21}. 
Moreover, the operator
$$
I+H(\psi_n)
$$
is invertible for each $n\ge 1$ and the inverses are uniformly bounded. Indeed, this can be seen by 
making a unitary transform,
$$
R_{\mu_n} H(\psi_n) R_{\mu_n}^*=  H(G_{\mu_n}\psi_n)= H(\hat{\psi}),
$$
using formula (\ref{RGH}). For the inverses we obtain the formula
$$
(I+H(\psi_n))\iv = R_{\mu_n}^* (I+H(\hat{\psi}))\iv R_{\mu_n},
$$
from which the uniform boundedness of the inverses follows.

(b):\
Let us start with some preliminary considerations. Observe that $P_n=H(t^n)^2$ and that 
$$
h_{R,n}:=G_{\mu_n}(t^n)= \left(\frac{t+\mu_n}{1+\mu_nt}\right)^n,
$$
which implies, by  (\ref{RGH}), that 
$$
R_{\mu_n} P_n R_{\mu_n}^*= H(h_{R,n})^2.
$$
The functions $h_{R,n}$ are uniformly bounded in the $L^\iy$-norm, and using the asymptotics
(\ref{mu.n}) it is easily seen that $h_{R,n}\to h_R$ in measure as $n\to\iy$. From Lemma \ref{l3.3}
we can conclude that 
$$
H(h_{R,n})\to H(h_R), \quad T(h_{R,n})\to T(h_R), \quad H(h_{R,n})^*\to H(h_R)^*
$$
strongly on $H^2(\T)$. Moreover, because the functions $h_{R,n}$ and all of its derivatives
converge locally uniformly on $\T\setminus\{-1\}$ to the function $h_{R}$ and its corresponding derivatives, we can conclude that 
$$
H(h_{R,n} f)\to H(h_R f)
$$
in the trace norm whenever $f$ is a sufficiently smooth function on $\T$ which vanishes identically in a neighborhood of $t=-1$. 

Furthermore, from the definitions it follows that 
$$
R_{\mu_n} H(\psi) R_{\mu_n}^*=H(\psi),\qquad
R_{\mu_n} H(\psi_n) R_{\mu_n}^*=H(\hat{\psi}).
$$
Combining all this we conclude that the first expression to analyze equals
$$
H( h_{R,n})^2 (I+H(\psi))\iv H( h_{R,n})^2 +I-H( h_{R,n})^2,
$$
and that this expression converges strongly to
$$
H(h_R)^2(I+H(\psi))\iv H(h_R)^2+I-H(h_R)^2.
$$
This is equal to the right hand side. By definition $H(h_R)^2= \Pi_R$.

The second expression to analyze equals
$$
H( h_{R,n})^2 \left((I+H(\hat{\psi}))\iv-(I+H(\psi))\iv
\right) H( h_{R,n})^2.
$$
Now choose two smooth, even functions $f_1$, $f_{-1}$ such that $f_1+f_{-1}=1$ and such that 
$f_1$ vanishes on a neighborhood of $1$, while $f_{-1}$ vanishes on a neighborhood of $-1$. 
Then the above expression equals the sum
$$
H( h_{R,n})^2 T(f_{-1})\left((I+H(\hat{\psi}))\iv-(I+H(\psi))\iv
\right)H( h_{R,n})^2\qquad
$$
$$
{} +
H( h_{R,n})^2 T(f_1)\left((I+H(\hat{\psi}))\iv-(I+H(\psi))\iv
\right)H( h_{R,n})^2.
$$
In view of the first term, we write (see (\ref{f.Hab}))
$$
H( h_{R,n}) T(f_{-1})= H( h_{R,n} f_{-1})-T( h_{R,n}) H(f_{-1}).
$$
This is trace class and converges in the trace norm to
$$
H(h_R) T(f_{-1})= H(h_R f_{-1})-T(h_R) H(f_{-1}).
$$
As to the second term, consider
$$
T(f_1)H(\psi)+H(f_1)T(\tilde{\psi})=H(f_1 \psi)=H(\psi f_1)=H(\psi)T(f_1)+T(\psi)H(f_1).
$$
Hence $T(f_1)H(\psi)$ equals $H(\psi)T(f_1)$ plus a trace class operator. From this we can conclude
that 
$$
T(f_1)(I+H(\psi))\iv=(I+H(\psi))\iv T(f_1)+\mbox{ trace class}.
$$
Next observe that $f_1(\psi-\hat{\psi})$ is continuous and has one-sided derivatives at $t=-1$, while it  is smooth elsewhere. 
Hence $H((\psi-\hat{\psi})f_1)$ is a trace class operator, which implies that 
$T(f_1)H(\psi)$ equals $T(f_1)H(\hat{\psi})$ plus a trace class operator. Using the above it follows that 
$H(\psi)T(f_1)$ equals $T(f_1)H(\hat{\psi})$ modulo trace class. This implies that 
$$
(I+H(\psi))\iv T(f_1)=T(f_1)(I+H(\hat{\psi}))\iv+\mbox{ trace class}.
$$
Combining with the above we see that 
$$
T(f_1)\left((I+H(\hat{\psi}))\iv-(I+H(\psi))\iv\right)
$$
is trace class, and from this the desired conclusion follows easily.
\end{proof}

The operator $\Pi_R$ is a projection, i.e., $\Pi_R^2=\Pi_R$. Indeed, we have
$h_R\in H^\iy(\T)$ and $\tilde{h}_R=h_R\iv\in\ovl{H^\iy(\T)}$. Hence by (\ref{f.Hab})
we have $H(h_R)T(h_R)=H(h_R \tilde{h}_R)-T(h_R)H(\tilde{h}_R)=0$, whence
by (\ref{f.Tab}),
$$
H(h_R)^3=H(h_R)(I-T(h_R)T(\tilde{h}_R))=H(h_R).
$$
This shows that $H(h_R)^2$ is a projection. Let us denote the image of
$\Pi_R$ by  $``\mathrm{im}\, \Pi_R$'',  which is a closed subspace of $H^2(\T)$.

\begin{theorem} \label{t6.2}
Let $R>0$ be fixed,  $|\Re \alpha|<1$,  and let $\psi,\psi_n,$ and $\hat{\psi}$ be defined by
(\ref{f.17}), (\ref{f.17b}), (\ref{psihat}), and (\ref{mu.n}). Then the operator
$$
\Big(\Pi_R(I+H(\psi))\iv\Pi_R\Big) \Big|_{\mathrm{im} \,\Pi_R}
$$
is invertible on the space $\mathrm{im} \, \Pi_R$, and 
\bqn
\lefteqn{
\lim_{n\to\iy}
\frac{\det\Big( P_n(I+H(\psi_n))\iv P_n\Big)}{\det\Big( P_n(I+H(\psi))\iv P_n\Big) }}
\hspace{4ex}\nn\\
&=&
\det\Big[\Big(\Pi_R(I+H(\psi))\iv\Pi_R\Big)\iv\Big(
\Pi_R(I+H(\hat{\psi}))\iv \Pi_R\Big)\Big|_{\mathrm{im}\, \Pi_R}\Big].\label{f.20}
\eqn
In particular, the expression on the right hand side is well-defined.
\end{theorem}
\begin{proof}
The expression for which we want to determine the limit is equal to the determinant of
$$
\Big(P_n(I+H(\psi))\iv P_n+Q_n\Big)\iv
\Big(P_n(I+H(\psi_n))\iv P_n+Q_n\Big).
$$
By Lemma \ref{l6.1}(a), the inverses of $I+H(\psi)$ and $I+H(\psi_n)$ exist for each $n$. Moreover,
by Corollary \ref{c3.8} the sequence $\{P_n(I+H(\psi))\iv P_n\}_{n\in\Z_+}$ is stable. In fact, we can even say that the inverses of these finite sections exist for each $n$ as noted in Theorem \ref{t4.3}. Hence the above expression makes sense for all $n$.

We rewrite the above expression as
$$
I+
\Big(P_n(I+H(\psi))\iv P_n+Q_n\Big)\iv
\Big(P_n\left( (I+H(\psi_n))\iv -(I+H(\psi))\iv\right) P_n\Big)
$$
and multiply with the unitaries $R_{\mu_n}^*$ and $R_{\mu_n}$, which does not change the 
value of the determinant, in order to obtain
$$
I+
\Big(R_{\mu_n}^* P_n (I+H(\psi))\iv P_n R_{\mu_n} +R_{\mu_n}^* Q_n R_{\mu_n} \Big)\iv
\Big(R_{\mu_n}^*P_n\left( (I+H(\psi_n))\iv -(I+H(\psi))\iv \right) P_nR_{\mu_n}\Big).
$$
By Lemma \ref{l6.1}(b) the expression on the right hand side tends to
$$
\Pi_R\left( (I+H(\hat{\psi}))\iv-(I+H(\psi))\iv\right) \Pi_R
$$
in the trace norm. Since the sequence $P_n (I+H(\psi))\iv P_n$ is stable,  so is the sequence
$$
R_{\mu_n}^* P_n (I+H(\psi))\iv P_n R_{\mu_n} +R_{\mu_n}^* Q_n R_{\mu_n}.
$$
This sequence converges strongly to
$$
\Pi_R(I+H(\psi))\iv \Pi_R+(I- \Pi_R).
$$
This last operator is also invertible. This follows either from Lemma  \ref{l.basic}(ii), or 
from Corollary \ref{c3.4}. Combining all this, we conclude that the limit under consideration exists and equals the determinant
of
$$
I+\Big(\Pi_R(I+H(\psi))\iv \Pi_R+(I- \Pi_R)\Big)\iv
\Big(\Pi_R\left( (I+H(\hat{\psi}))\iv-(I+H(\psi))\iv\right) \Pi_R
\Big).
$$
This operator equals
$$
\Big(\Pi_R(I+H(\psi))\iv \Pi_R+(I- \Pi_R)\Big)\iv
\Big(\Pi_R(I+H(\hat{\psi}))\iv  \Pi_R+(I- \Pi_R)\Big),
$$
the determinant of which is the same as in (\ref{f.20}).
\end{proof}

We remark that the determinant on the right hand side of (\ref{f.20}) cannot be written as 
the quotient of two determinants of the corresponding operators because (as one can show) 
neither of the operators
$$
\Big(\Pi_R(I+H(\psi))\iv\Pi_R\Big)\quad \mbox{ and }\quad 
\Big(\Pi_R(I+H(\hat{\psi}))\iv \Pi_R\Big)
$$
is of the form identity plus trace class. The reason is that the generating functions of the Hankel operators have a jump discontinuity at $t=-1$.

Let us now state the main result of this and the previous section. It concludes  Step 1
discussed in the introduction.

\begin{corollary}\label{c6.3}
Let $R>0$ and $|\Re\alpha|<1$. Then 
\be\label{f.58}
P^{(\alpha)}(R)= \exp\Big(-\frac{R^2}{4}+\alpha R\Big)
\det\Big[\Big( P_R(I+H_\R(\psi))\iv P_R\Big)\iv\Big(
 P_R(I+H_\R(\hat{\psi}))\iv  P_R\Big)\Big],
\ee
where
$$
\psi(x)=\uh_{-1/2-\alpha,0}(x)\uh_{1/2+\alpha,\iy}(x),\qquad
\hat{\psi}(x) =\uh_{-1/2,0}(x)\uh_{1/2+\alpha,\iy}(x).
$$
\end{corollary}
\begin{proof}
By Corollary \ref{c22}, the operators $I+H_\R(\psi)$ and
$I+H_\R(\hat{\psi})$ are invertible on $L^2(\R_+)$. Moreover, the operator
$P_R(I+H_\R(\psi))\iv P_R$ is invertible on $L^2[0,R]$
by Corollary \ref{c3.4}.
In the previous theorem we have shown that 
$$
A=\Big[\Big(\Pi_R(I+H(\psi))\iv\Pi_R\Big)\iv\Big(
\Pi_R(I+H(\hat{\psi}))\iv \Pi_R\Big)\Big|_{\mathrm{im}\, \Pi_R}\Big].
$$
(with $\psi$ and $\hat{\psi}$ as defined by (\ref{f.17}) and (\ref{psihat})) makes sense and is of the form identity  plus trace class.
We apply the transform $\cS$ to $A$ and obtain that 
$$
\cS A \cS^{-1} = 
\Big( P_R(I+H_\R(\psi))\iv P_R\Big)\iv\Big(
 P_R(I+H_\R(\hat{\psi}))\iv  P_R\Big),
$$
which thus is also of the form identity plus trace class. Hence the determinant 
in (\ref{f.58}) is well defined and equals the right hand side in (\ref{f.20}).
Notice that here we have used 
$$
\cS H(a) \cS^{-1}=H_\R(\ah),\quad \cS \Pi_R \cS\iv = P_R,
$$
see (\ref{HHR}) and (\ref{uuh}), where the last identity follows from the fact that 
$P_R= H_\R(e^{ixR})^2$, and $\cS H(h_R)\cS^{-1}= H_\R(e^{ixR})$.

Combining this with the previous theorem and Theorem \ref{t4.3}, it follows that 
the identity (\ref{f.58}) holds for all real $\alpha$ with $|\alpha|<1$.

Now, by Theorem \ref{t4.1} the quantity $P^{(\alpha)}(R)$ depends analytically on $\alpha$
for $|\Re \alpha|<1$. The same is true for the right hand side of (\ref{f.58}) since the generating functions
$\psi$ and $\hat{\psi}$ depend analytically on $\alpha$. Hence the identity 
(\ref{f.58}) holds not just for real $\alpha$, but for complex $\alpha$.
This completes the proof.
\end{proof}

%%%%%%%%%%%%%%%%%%%%%%%%%%%%%%%%%%%%%%%%%%%%%%%%%%%%%%%%%%%%%%%%%%%%%%%

\section{Localization}
\label{s7}

Let $\cB$ stand for the Banach algebra of all function $a\in L^\iy(\T)$ for which both 
$H(a)$ and $H(\ta)$ are trace class (see \cite[Sec.~10.2-3]{BS}). The norm in $\cB$  is defined as 
$$
\|a\|_{\cB}=|a_0|+\|H(a)\|_{\cC_1(H^2(\T))}+\|H(\ta)\|_{\cC_1(H^2(\T))}.
$$
The class $\cB$ can be identified with the so-called Besov class $B^1_1$, however, we will not make use of this fact. What is important for us is that $\cB$ contains all smooth functions as a dense subset and that the Riesz projection is bounded on $\cB$. Using this and Gelfand theory, it follows that the maximal ideal space of $\cB$ can be identified with $\T$. In other word, $a\in\cB$ is invertible in $\cB$ if and only if $a(t)\neq0$ for all $t\in\T$. Moreover, $a\in\cB$ possesses a logarithm in $\cB$ if and only if
$a$ possesses a continuous logarithm. 
Now we define the unital Banach algebras
$$
\cB_+=\cB\cap H^\iy(\T),\qquad \cB_-=\cB\cap \overline{H^\iy(\T)},
$$
and we can introduce the notion of Wiener-Hopf factorization in $\cB$ similar as for $\W$
(see Sec.~\ref{s4}). We also recall that $GB$ stands for the group of all invertible elements in a 
unital Banach algebra $B$.

What is also important to us is that $\cB$ contains all continuous functions on $\T$
which are smooth except at finitely many points at which the one-sided derivatives exist. This is easy to prove using the fact that the Fourier coefficients of such a function decay as $O(n^{-2})$.

We start with the following result, which has essentially been proven already in 
\cite{BE1}. We will sketch the main idea of the proof in order to indicate how the constant arises. 
For more details we refer to \cite{BE1}.

\begin{proposition}\label{p5.1}
Let $b_+\in G\cB_+$. Then 
\be\label{f.29}
\det(I+H(b_+\tb_+\iv))=
\left(\frac{b_+(1))}{b_+(-1)}\right)^{1/2}\exp\Big(-\frac{1}{2}\sum_{k=1}^\iy k [\log b_+]_k^2\Big).
\ee
\end{proposition}
\begin{proof} 
Let $b=b_+\tb_+$. Then $T(b)=T(\tb_+)T(b_+)$, whence $T\iv(b)=T(b_+\iv)T(\tb_+\iv)$, and
$$
T\iv(b) H(b)=T(b_+\iv)T(\tb_+\iv)H(b_+ \tb_+)=T(b_+\iv)H(b_+).
$$
Using $H(b_+)T(b_+\iv)=H(b_+\tb_+\iv)$, it follows that 
$$
\det (I+T\iv(b)H(b))=\det(I+H(b_+\tb_+\iv)).
$$
Here we have used the formulas (\ref{f.THabc}).
In \cite[Thm.~2.5]{BE1} it was shown that $\det (I+T\iv(b)H(b))$ equals the 
right hand side in (\ref{f.29}). The crucial point 
is to introduce an operator $M(a)=T(a)+H(a)$ and to observe that 
$M(a_1a_2)=M(a_1)M(a_2)$ whenever $a_2=\ta_2$. Indeed, this last identity follows from (\ref{f.Tab}) and (\ref{f.Hab}). 

Now consider the function $b_\lambda=\exp(\lambda \log b)$, which depends
analytically on $\lambda$, and define the analytic function $f(\lambda)=\det T\iv(b_\lambda) M(b_\lambda)$.
Take the logarithmic derivative of $f$ (by employing the formula
$(\log\det F)'=\tr(F'F\iv)$) and differentiate once more. A simple computation yields
$$
\frac{f'(\lambda)}{f(\lambda)}= \tr\left(M(b_\lambda')M\iv(b_\lambda)-T\iv(b_\lambda)T(b_\lambda')\right)
=\tr\left(M(\log b)-T\iv(b_\lambda)T(b_\lambda')\right)
$$
and
$$
\left(\frac{f'(\lambda)}{f(\lambda)}\right)'=\tr\left(  T\iv(b_\lambda)T(b_\lambda') T\iv(b_\lambda)T(b_\lambda') -T\iv(b_\lambda)T(b_\lambda'')\right).
$$
Now we use the facts that $b_\lambda'= b_\lambda \log b $ and $b_\lambda''= b_\lambda \log^2b$
and that $b_\lambda$ has a canonical Wiener-Hopf factorization in $\cB$, say $b_\lambda=b_{\lambda,+} b_{\lambda,-}$ with $b_{\lambda,\pm}\in G\cB_\pm$. Using $T\iv(b_\lambda)=T(b_{\lambda,+}\iv)T(b_{\lambda,-}\iv)$
and (\ref{f.THabc}) we obtain 
$$
\left(\frac{f'(\lambda)}{f(\lambda)}\right)'=\tr\left(  T(\log b) T(\log b)-T(\log ^2 b)\right)=
-\tr \left(H(\log b) H(\log b)\right),
$$
which does not depend on $\lambda$. Integrating and fixing the constants with the values of $f$ and $f'$ at $\lambda=0$ yields
$$
f(\lambda)=\exp\left(\lambda\,  \tr \,H(\log b)-\frac{\lambda^2}{2}\tr  \left(H(\log b) H(\log b)\right)\right). 
$$
{}From this  the assertion follows by putting $\lambda=1$ and by evaluating the traces.
\end{proof}

We proceed now with two rather technical lemmas. Let us recall the notation 
$G_r$ defined in (\ref{f.Gm}).

\begin{lemma}\label{l7.2a}
Let $a=u_{\alpha,1}$, $c=u_{\gamma,1}-1$, and  $\psi_r(t)=(1-t)a(t)(G_rc)(t)$ for $t\in\T$, $r\in[0,1)$.
Then $H(\psi_r)$ is trace class and tends to zero in the trace norm as $r\to1$.
\end{lemma}
\begin{proof}
We first remark that it is rather easy to see that $\psi_r$ converges uniformly to zero on $\T$.
In fact, the function $G_r c$ converges to zero locally uniformly on $\T\setminus\{1\}$ and it is  uniformly bounded on $\T$. Hence the multiplication with $(1-t)$ implies the uniform convergence on all of $\T$. As a consequence the Hankel operator converges in the operator norm to zero. The technical difficulty is
the convergence in the trace norm.

For fixed $r$,  the functions $\psi_r$ are smooth on $\T\setminus\{1\}$, continuous on $\T$, and have one-sided higher order derivatives at $t=1$. This implies that all the Hankel operators
are trace class. We will soon use the fact that the functions $(1-t)a(t)$ and $(1-t) c(t)$ have the same properties. Hence the Hankel operators with these symbols are trace class, too.

In order to show the convergence of $H(\psi_r)$ in the trace norm use (\ref{f.Hab}) and write
$$
H(\psi_r) = T((1-t)a)H(G_rc)+H((1-t)a)T(\wt{G_r c}).
$$
The sequence $G_rc$ converges in measure to zero and is uniformly bounded. Hence the adjoint of
$T(\wt{G_r c})$ converges strongly to zero. On the other hand, the Hankel operator
$H((1-t)a)$ is trace class. It follows that the second term in the above expression converges to zero in the trace norm. Hence our concern from now on is the first term. 

Applying the unitary operator $R_r$ defined in (\ref{f.Rm}), the trace norm of this operator is equal to the trace norm of
\bqn
R_r^* T((1-t)a)H(G_rc) R_r &=&
R_r^* T((1-t)a) R_r R_r^* H(G_rc) R_r \nn\\
&=& T(G_r\iv((1-t)a))H(c).\nn
\eqn
One computes easily that 
$$
G_r\iv:(1-t)\mapsto \frac{(1-r)(1-t)}{1+rt}.
$$
Because of (\ref{f.THabc}) we conclude that 
\bqn\label{fx.53}
T(G_r\iv((1-t)a))H(c) = T\left(\frac{1-r}{1+rt}\,G_r\iv a\right)T(1-t)H(c).
\eqn
Observe that $|(1-r)/(1+rt)|\le1$ and that $(1-r)/(1+rt)$ converges locally uniformly on $\T\setminus\{-1\}$ to zero. Hence this term converges to zero in measure. The sequence $G_r\iv a$ is also uniformly bounded. We conclude that the Toeplitz operator on the left of the right hand side of the last equation tends strongly to zero. Using (\ref{f.Hab}) now write
$$
T(1-t)H(c) = H((1-t)c)-H(1-t)T(\tilde{c})
$$ 
to see that this operator is trace class. We conclude that the expression (\ref{fx.53}) converges to zero in the trace norm. This finishes the proof.
\end{proof}

Let us define the expression
\bqn\label{f.54}
K(a,b)= (I+H(ab))-(I+H(a))(I+H(b)).
\eqn
for $a,b\in L^\iy(\T)$.

\begin{lemma}\label{l5.2}
Let $a=u_{\alpha,1}$, $b=u_{\beta,-1}$, and
$$
a_r(t)=\left(\frac{1-rt}{1-rt\iv}\right)^{\alpha},\qquad
b_r(t)=\left(\frac{1+rt}{1+rt\iv}\right)^{\beta},\qquad r\in [0,1).
$$
Then $K(a,b)$ is trace class, and  $K(a_r,b_r)\to K(a,b)$ in the trace norm as $r\to 1$.
\end{lemma}
\begin{proof}
We first observe that we can write $K(a,b)$ as
$$%\be\label{f.24}
K(a,b)=H((a-1)(b-1))-H(a)H(b).
$$%\ee

The function $(a-1)(b-1)$ is continuous and sufficiently smooth on $\T\setminus\{\pm1\}$ and has one-sided derivatives of arbitrary high order at $t=\pm1$. Hence the Hankel operator $H((a-1)(b-1))$ is trace class. To see that $H(a)H(b)$ is trace class, decompose $1=f_1+f_{-1}$ into smooth and even functions such that $f_1$ is identically zero in a neighborhood of $1$, while $f_{-1}$ is identically zero
in  a neighborhood of $-1$. Then write, using (\ref{f.Hab}),
\bqn
H(a)H(b) &=& H(a)T(f_1)H(b)+H(a)T(f_{-1})H(b)\nn\\
&=&
\Big(H(a f_1)-T(a)H(f_1)\Big)H(b) + H(a) \Big( H(f_{-1}b)-H(f_{-1})T(\tb)\Big),
\label{f.25} 
\eqn
which is trace class because so are the Hankel operators with the symbols
$af_1$, $f_{-1}b$, $f_1$, $f_{-1}$.

The functions $a_r$ and $b_r$ are smooth. Hence the Hankel operators
$H((a_r-1)(b_r-1))$ and $H(a_r)H(b_r)$ are trace class for each $r$. Let us first show that 
$$
H(a_r)H(b_r)\to H(u_{\alpha,1})H(u_{\beta,-1}).
$$
We use the identity (\ref{f.25}) with $a_r$ and $b_r$ instead of $a$ and $b$, and hence we
have
$$
H(a_r)H(b_r) = \Big(H(a_r f_1)-T(a_r)H(f_1)\Big)H(b_r) + H(a_r) \left( H(f_{-1}b_r)-H(f_{-1})T(\tb_r)\right).
$$
Now remark that $a_r$ and $b_r$ are uniformly bounded and converge in measure to $a$ and $b$.
Hence we have following strong convergences,
 $$
 T(a_r)\to T(a),\quad H(b_r)^*\to H(b)^*,\quad H(a_r)\to H(a),\quad T(\tb_r)^*\to T(\tb)^*.
 $$
On the other hand, $a_r f_1\to a f_1$ and $b_r f_{-1}\to b f_{-1}$, e.g., in the norm of
$C^2(\T)$, whence it follows that $H(a_r f_1)\to H(a f_1)$ and $H(f_{-1} b_r)\to H(f_{-1}b)$
in the trace norm. Combining all this we can conclude that $H(a_r)H(b_r)\to H(a)H(b)$ in the trace norm.

In order to treat convergence of the Hankel operator $H((a_r-1)(b_r-1))$ we use the above functions
$f_{1}$ and $f_{-1}$, and decompose
$$
H((a_r-1)(b_r-1))=H((a_r-1)f_{1}(b_r-1))+H((a_r-1)f_{-1}(b_r-1)).
$$
Without loss of generality it suffices to consider the last term. In fact, the first term on the right can be 
transformed into the same kind of expression by the unitary operator
$Y_{-1}:f(t)\mapsto f(-t)$, $t\in\T$. Now write, using (\ref{f.Hab}),
$$
H((a_r-1)f_{-1}(b_r-1))=H((a_r-1)(t-1))T\Big(f_{-1}\frac{\tilde{b}_r-1}{t\iv-1}\Big)+
T((a_r-1)(t-1))H\Big(f_{-1}\frac{b_r-1}{t-1}\Big)
$$

Let us first focus on the terms containing $b_r$. Let $K$ be any compact subset of $\T\setminus\{-1\}$.
We may think of the functions $b_r(t)$ and $b(t)$ as being defined and analytic on a suitable
open neighborhood of $K\subset \C$. 
On this neighborhood of $K$,
we have uniform convergence $b_r\to b$. Because $b_r(1)=b(1)=1$, we have also uniform convergence
$$
\frac{b_r(t)-1}{t-1}\to \frac{b(t)-1}{t-1}
$$
along with all derivatives. It follows that 
$$
T\Big(f_{-1}\frac{\tilde{b}_r-1}{t\iv-1}\Big)\to T\Big(f_{-1}\frac{\tilde{b}-1}{t\iv-1}\Big),\qquad
H\Big(f_{-1}\frac{b_r-1}{t-1}\Big)\to H\Big(f_{-1}\frac{b-1}{t-1}\Big).
$$
where the convergence of the Toeplitz operators is in the operator norm and the convergence of the Hankel operators is in the trace norm. 

Next observe that $a_r\to a$ locally uniformly on $\T\setminus\{1\}$, whence it follows that 
$(t-1)a_r(t)\to (t-1)a(t)$ uniformly on $\T$. Hence
$$
T((a_r-1)(1-t))\to T((a-1)(1-t))
$$
in the operator norm. In order to treat the Hankel operators with the same symbols write
$$
a_r(t)= \left(\frac{1-rt}{1-rt\iv}\right)^\alpha=(-t)^{\alpha}\left(-\frac{t-r}{1-rt}\right)^{-\alpha}=
u_{\alpha,1}(t)(G_r u_{-\alpha,1})(t).
$$
Now Lemma \ref{l7.2a} implies that $H((t-1)a_r)\to H((t-1)a)$ in the trace norm. 
 Combining all the previous considerations, it follows that 
$$
H((a_r-1)f_{-1}(b_r-1))\to H((a-1)f_{-1}(b-1))
$$
in the trace norm, which completes the proof.
\end{proof}

\begin{proposition}\label{p7.3}
Let $a=u_{-1/2-\alpha,1}$, $b=u_{1/2+\beta,-1}$, and assume  $|\Re\alpha|<1$ and $|\Re\beta|<1$.
Then the operator determinant 
\be\label{f.d26}
\det\Big((I+H(a))\iv (I+H(ab)) (I+H(b))\iv \Big)
\ee
is well defined and equals  $2^{-(1/2+\alpha)(1/2+\beta)}$.
\end{proposition}
\begin{proof}
Theorem \ref{t21} implies that the inverses of $I+H(a)$, $I+H(b)$, and $I+H(ab)$ exist.
Defining $K(a,b)$ as in (\ref{f.54}) we can write
$$
(I+H(a))\iv (I+H(ab))(I+H(b))\iv = I+(I+H(a))\iv K(a,b) (I+H(b))\iv,
$$
and hence this operator is of the form identity plus trace class. Its determinant is well defined.

In order to compute the value of this operator determinant, we 
approximate the functions $a$ and $b$ by smooth functions
$$
a_r(t)=\left(\frac{1-rt}{1-rt\iv}\right)^{-1/2-\alpha},\qquad
b_r(t)=\left(\frac{1+rt}{1+rt\iv}\right)^{1/2+\beta},\qquad r\in [0,1),
$$
and then let $r\to 1$. 

The sequences $a_r$ and $b_r$ are bounded in the $L^\iy$-norm and converge in measure to the 
functions $a$ and $b$, respectively. This implies that 
$H(a_r)\to H(a)$ and $H(b_r)\to H(b)$ strongly on $H^2(\T)$ as $r\to1$, and the same holds for the adjoints. Theorem \ref{t.23st} implies that the sequences $\{I+H(a_r)\}_{r\in[0,1)}$ and $\{I+H(b_r)\}_{r\in[0,1)}$ are stable. Hence
using Lemma \ref{l.basic}(iii) it follows that 
$$
(I+H(a_r))\iv \to (I+H(a))\iv \quad\mbox{and}\quad 
((I+H(b_r))\iv)^*\to ((I+H(b))\iv)^*
$$ 
strongly on $H^2(\T)$. This together with Lemma \ref{l5.2} implies
that
$$
(I+H(a_r))\iv K(a_r,b_r) (I+H(b_r))\iv\to (I+H(a))\iv K(a,b) (I+H(b))\iv
$$
in the trace norm as $r\to 1$. Hence the determinant
\bqn\label{f.50}
\det\Big((I+H(a_r))\iv (I+H(a_rb_r)) (I+H(b_r))\iv \Big)
\eqn
converges to  (\ref{f.d26}) as $r\to 1$. In (\ref{f.50}) each of the Hankel operators is trace class,
and hence we can split it into the product/quotient of three determinants, each of which we can evaluate
by Proposition \ref{p5.1}. We obtain
$$
\det(I+H(a_r))=\left(\frac{a_{+,r}(1)}{a_{+,r}(-1)}\right)^{1/2}
\exp\Big(-\frac{1}{2}\sum_{k=1}^\iy k [\log a_{+,r}]_k^2\Big),
$$
$$
\det(I+H(b_r))=
\left(\frac{b_{+,r}(1)}{b_{+,r}(-1)}\right)^{1/2}\exp\Big(-\frac{1}{2}\sum_{k=1}^\iy k [\log b_{+,r}]_k^2\Big),
$$
$$
\det(I+H(a_{r}b_r))=
\left(\frac{a_{+,r}(1)b_{+,r}(1)}{a_{+,r}(-1)b_{+,r}(-1)}\right)^{1/2}\exp\Big(-\frac{1}{2}\sum_{k=1}^\iy k [\log a_{+,r}+\log b_{+,r}]_k^2\Big).
$$
Here $a_{+,r}(t)=(1-rt)^{-1/2-\alpha}$, $b_{+,r}(t)=(1+rt)^{1/2+\beta}$.
Hence the determinant (\ref{f.50}) equals
\bqn
\exp\left(-\sum_{k=1}^\iy k [\log a_{+,r}]_k [\log b_{+,r}]_k\right) &=&
\exp\left(\gamma \sum_{k=1}^\iy k \left( -\frac{r^{k}}{k}\right) \left(-\frac{(-r)^{k}}{k}
\right)\right)\nn\\
&=&
\exp\left(\gamma \sum_{k=1}^\iy \frac{(-r^2)^k}{k}\right) =(1+r^2)^{-\gamma}\nn
\eqn
with $\gamma=(1/2+\alpha)(1/2+\beta)$.
Now take $r\to1$.
\end{proof}

\begin{proposition}\label{p7.4}
Let $a=\uh_{-1/2-\alpha,0}$, $b=\uh_{1/2+\beta,\iy}$, and assume that 
$|\Re\alpha|<1$, $|\Re\beta|<1$. Then we can write
$$
 P_R(I+H_\R(ab))\iv  P_R=  P_R(I+H_\R(b))\iv  P_R(I+H_\R(a))\iv  P_R+  P_R K  P_R+C_R,
$$
where $K$ is a trace class operators on $L^2(\R_+)$ and $C_R$ are  trace class operators on $L^2[0,R]$ tending to zero in the trace norm as $R\to\iy$.
\end{proposition}
\begin{proof}
Let us first remark that the operators $I+H_\R(a)$, $I+H_\R(b)$, and $I+H_\R(ab)$ are invertible because of Corollary \ref{c22}. From the part of Lemma \ref{l5.2} which states that 
$K(a,b)$ is trace class (see also (\ref{f.54})) it follows that
\be\label{f.K}
(I+H_\R(ab))\iv = (I+H_\R(b))\iv (I+H_\R(a))\iv +K,
\ee
where $K$ is a trace class operator. To see this
we have to apply the transformation $\cS$ (see (\ref{HHR}) and (\ref{uuh})).

Let $V_{\pm R}=W(e^{\pm iRx})$. These operators are the forward and backward shifts,
$$
(V_R f)(x)=\left\{\ba{cl} f(x-R) &\mbox{ if } x>R\\ 0 & \mbox{ if } 0\le x\le R,\ea\right.
\qquad 
(V_{-R}f)(x)= f(x+R),\quad x\ge0.
$$
Clearly, $Q_R=V_RV_{-R}$, and the formula $V_{-R}H_\R(c)=H_\R(e^{-iRx}c)=H_\R(c)V_R$ holds,
which follows from the continuous analogue of (\ref{f.THabc}). Using the identities
$$
(I+B)\iv= I-(I+B)\iv B,\qquad
(I+A)\iv =I-A(I+A)\iv
$$
with $B=H_\R(b)$ and $A=H_\R(a)$, it follows that
$$
P_R (I+H_\R(b))\iv Q_R (I+H_\R(a))\iv P_R
= P_R (I+H_\R(b))\iv H_\R(b) Q_R H_\R(a)(I+H_\R(a))\iv P_R
$$
$$
= P_R (I+H_\R(b))\iv V_{-R} H_\R(b)  H_\R(a) V_R (I+H_\R(a))\iv P_R.
$$
This term converges to zero because $H_\R(a)H_\R(b)$ is trace class and $V_{-R}\to 0$ strongly
as $R\to \iy$. Combining the previous formulas, using $I=P_R+Q_R$, the desired formula follows easily.
\end{proof}

\begin{theorem}\label{t7.5}
Let $a=\uh_{-1/2-\alpha,0}$, $b=\uh_{1/2+\beta,\iy}$, and assume that 
$|\Re\alpha|<1$, $|\Re\beta|<1$. Then
$$
\lim_{R\to\iy}
\frac{\det\Big[ \Big(P_R(I+H_\R(b))\iv P_R\Big)\iv  \Big(P_R(I+H_\R(ab))\iv P_R\Big)\Big] }{
\det\Big(P_R(I+H_\R(a))\iv P_R\Big)}=2^{(1/2+\alpha)(1/2+\beta)}
$$
All the expressions on the left hand side are well defined for sufficiently large $R$.
\end{theorem}
\begin{proof}
Because of Corollary \ref{c3.7} the sequences
$$
\{P_R(I+H_\R(a))\iv P_R\}_{R>0}\quad\mbox{and}\quad 
\{P_R(I+H_\R(b))\iv P_R\}_{R>0}
$$
are stable. Moreover, taking the strong limit on $L^2(\R_+)$ as $R\to\iy$ we get $(I+H_\R(a))\iv$
and $(I+H_\R(b))\iv$. Similarly, we can take the strong limit of the adjoints.

For sufficiently large $R$ it is thus possible to consider
$$
\Big(P_R(I+H_\R(b))\iv P_R\Big)\iv
\Big(P_R(I+H_\R(ab))\iv P_R\Big)\Big(P_R(I+H_\R(a))\iv P_R\Big)\iv ,
$$
which, by Proposition \ref{p7.4}, equals
$$
P_R+\Big(P_R(I+H(b))\iv P_R\Big)\iv K\Big(P_R(I+H(a))\iv P_R\Big)\iv + \wt{C}_R
$$
with $\{\wt{C}_R\}_{R>0}$ being a sequence of trace class operators on $L^2[0,R]$ tending to zero in the trace norm, and with $K$ being trace class operator on $L^2(\R_+)$.
Complementing with the projection  $Q_R$,
which does not change the value of the determinant,  
we can rewrite this identity as
$$
Q_R+\Big(P_R(I+H_\R(b))\iv P_R\Big)\iv
\Big(P_R(I+H_\R(ab))\iv P_R\Big)\Big(P_R(I+H_\R(a))\iv P_R\Big)\iv ,
$$
$$
=I+(I+H_\R(b))K(I+H_\R(a))+\wt{D}_R.
$$
with a sequence $\{\wt{D}_R\}_{R>0}$ of trace class operators on $L^2(\R_+)$ tending to zero in the trace norm.
Using the expression of $K$ given in (\ref{f.K}), which can also be obtained by passing to the limit
$R\to\iy$ in the previous equation, it follows that the above equals
$$
(I+H_\R(b))(I+H_\R(ab))\iv(I+H_\R(a))+\wt{D}_R.
$$
It follows that the following operator determinant is well defined
$$
\det\Big[ \Big(P_R(I+H_\R(b))\iv P_R\Big)\iv
\Big(P_R(I+H_\R(ab))\iv P_R\Big)\Big(P_R(I+H_\R(a))\iv P_R\Big)\iv\Big] ,
$$
and that its limit $R\to\iy$ equals
$$
\det\Big( (I+H_\R(b))(I+H_\R(ab))\iv(I+H_\R(a))\Big)
$$
Applying the transform $\cS$ (see (\ref{HHR}) and (\ref{uuh})) and noting that the reciprocal of resulting determinant has been computed in Proposition \ref{p7.3}, this completes the proof.
\end{proof}

%%%%%%%%%%%%%%%%%%%%%%%%%%%%%%%%%%%%%%%%%%%%%%%%%%%%%%%%%%%%%%%%%%%%%%%

\section{The final result}
\label{s8}

Let us now put all the pieces together and derive the final result.
In Corollary \ref{c6.3} we have shown that for $|\Re\alpha|<1$ and $R>0$ we have the identity
\bqn
&& P^{(\alpha)}(R) =
\exp\Big(-\frac{R^2}{4}+\alpha R\Big)\times
\nn\\
&&\qquad
\det\Big[\Big( P_R(I+H_\R(\uh_{-1/2-\alpha,0}\uh_{1/2+\alpha,\iy}))\iv P_R\Big)\iv\Big(
 P_R(I+H_\R(\uh_{-1/2,0}\uh_{1/2+\alpha,\iy}))\iv  P_R\Big)\Big].\nn
\qquad
\eqn
Therein the operator determinant is well defined (see also Corollary \ref{c3.4}).
For sufficiently large $R$, this determinant can be written as a product of the following two determinants,
\bqn
d_1(R) &=& \det\Big[\Big(P_R(I+H_\R(\uh_{-1/2-\alpha,0}\uh_{1/2+\alpha,\iy}))\iv P_R\Big)\iv\Big(
P_R(I+H_\R(\uh_{1/2+\alpha,\iy}))\iv P_R\Big)\Big] ,\nn
\\
d_2(R) &=&
\det\Big[\Big(P_R(I+H_\R(\uh_{1/2+\alpha,\iy}))\iv P_R\Big)\iv\Big(
P_R(I+H_\R(\uh_{-1/2,0}\uh_{1/2+\alpha,\iy}))\iv P_R\Big)\Big] .\nn
\eqn
In this connection observe that Corollary \ref{c3.7}  guarantees that the inverses of the various operators
$P_R(I+H_\R(\ast))\iv P_R|_{L^2[0,R]}$ exist. In particular, $d_k(R)\neq0$. Notice also that Theorem \ref{t7.5}
makes sure that both operator determinants are well defined. Moreover, Theorem \ref{t7.5} 
implies that 
$$ 
d_1(R) d_2(R) \sim 2^{-(1/2+\alpha)^2} 2^{(1/2+\alpha)/2} 
\frac{\det\Big( P_R(I+H_\R(\uh_{-1/2,0}))\iv P_R\Big)}
{\det \Big(P_R(I+H_\R(\uh_{-1/2-\alpha,0}))\iv P_R\Big)}
$$
as $R\to\iy$.
To complete the argument we need a (non-trivial) result established by E.L. Basor and the author in \cite[Sec.~3.6]{BE3}. Therein $G(z)$ stands for the Barnes $G$-function (\ref{Barnes}) (see also \cite{Bar}).

\begin{theorem}\label{t8.1}
Let $-3/2<\Re\gamma<1/2$. Then, as $R\to\iy$,
\bqn
\det \Big (P_R (I+H_\R(\hat{u}_{\gamma,0}))\iv P_R \Big)
\sim R^{\gamma^2/2+\gamma/2}(2\pi)^{-\gamma/2}2^{-\gamma^2-\gamma/2}\frac{G(1/2)}{G(1/2-\gamma)}.
\eqn
\end{theorem}

Applying  this theorem we obtain
$$ 
d_1(R) d_2(R) \sim 2^{-\alpha(1/2+\alpha)}\cdot
\frac{R^{-1/8} (2\pi)^{1/4}}
{R^{(\alpha^2-1/4)/2}(2\pi)^{\alpha/2+1/4} 2^{-\alpha(\alpha+1/2)} }
\cdot\frac{G(1+\alpha)}{G(1)}.
$$
Thus, after simplifying, we get our final result, which confirms the conjecture of Tracy and Widom \cite{TW94b}.

\begin{theorem}\label{t8.2}
Let $|\Re\alpha|<1$. Then, as $R\to\iy$,
\bqn
P^{(\alpha)}(R) \sim \exp\Big( -\frac{R^2}{4}+\alpha R-\frac{\alpha^2}{2}\log R\Big)
\frac{G(1+\alpha)}{(2\pi)^{\alpha/2}}.
\eqn
\end{theorem}

%%%%%%%%%%%%%%%%%%%%%%%%%%%%%%%%%%%%%%%%%%%%%%%%%%%%%%%%%%%%%%%%%%%%%%%

\end{document}